\pdfoutput=1
\documentclass[
final
]{dmtcs-episciences}

\usepackage[utf8]{inputenc}
\usepackage{subfigure}
\usepackage{amssymb}
\usepackage{amsmath,amsfonts,amsthm}
\usepackage{hyperref}
\usepackage{color}
\usepackage{cite}
\usepackage{tikz}
\definecolor{uuuuuu}{rgb}{0.26666666666666666,0.26666666666666666,0.26666666666666666}

\newtheorem{theorem}{Theorem}[section]

\newtheorem{lemma}[theorem]{Lemma}

\newtheorem{proposition}[theorem]{Proposition}

\usepackage[round]{natbib}

\author[M. Gholami and Y. Rowshan]{Mostafa Gholami\affiliationmark{1}\thanks{Corresponding Author.}
  \and Yaser Rowshan\affiliationmark{1}}
\title[The bipartite Ramsey numbers $BR(C_8, C_{2n})$]{THE BIPARTITE RAMSEY NUMBERS $BR(C_8, C_{2n})$}
\affiliation{
  Department of Mathematics, Institute for Advanced Studies in Basic Sciences (IASBS), Zanjan 66731-45137, Iran}
\begin{document}
\publicationdata
{vol. 25:2 }
{2023}
{15}
{10.46298/dmtcs.11207}
{2023-04-18; 2023-04-18; 2023-08-17; 2023-09-12}
{2023-09-14}
\maketitle
\begin{abstract}
 For the given bipartite graphs $G_1,G_2,\ldots,G_t$, the multicolor bipartite Ramsey number $BR(G_1,G_2,\ldots,G_t)$ is the smallest positive integer $b$, such that any $t$-edge-coloring of $K_{b,b}$ contains a monochromatic subgraph isomorphic to $G_i$, colored with the $i$th color for some $1\leq i\leq t$. We compute the exact values of the bipartite Ramsey numbers $BR(C_8,C_{2n})$ for $n\geq2$.
\end{abstract}

\section{Introduction}
Since 1956, when Erdös and Rado published the fundamental paper \cite{erdos1956partition}, Ramsey theory has grown into one of the most active areas of research in combinatorics while interacting with graph theory, number theory, geometry, and logic \cite{rosta2004ramsey}. Ramsey theory has many applications in several branches of mathematics. We refer to \cite{graham1990ramsey,R} to see these diverse applications. In particular, one can see the bipartite Ramsey numbers have many applications, and this motivated us to conduct a study on bipartite Ramsey numbers. The classical Ramsey number for the given numbers $n_1,\ldots,n_k$, is the smallest integer $n$ in a way that there is some $1\leq i\leq k$ for each  $k$-coloring of the edges of complete graph $K_n$, such that there is a complete subgraph of size $n_i$ whose edges are all the $i$th color.  However, there is no loss to work in any class of graphs and their subgraphs instead of complete graphs. If $G_1,G_2,\ldots,G_t$ be bipartite graphs, the multicolor bipartite Ramsey number $BR(G_1,G_2,\ldots,G_t)$ is defined as the smallest positive integer $b$, such that any $t$-edge-coloring of the complete bipartite graph $K_{b,b}$ contains a monochromatic subgraph isomorphic to $G_i$, colored with the $i$th color for some $i$. The existence of such a positive integer is guaranteed by a result of Erd\H{o}s and Rado \cite{erdos1956partition}. Recently, new versions of Ramsey numbers, such as multipartite Ramsey numbers have been defined. One can refer to \cite{day2001multipartite}, \cite{rowshan2021size}, \cite{rowshan2023multicolor}, \cite{rowshanbipartite} and their references for further studies.\\

The exact values of the bipartite Ramsey numbers $BR(P_n,P_m)$ of two paths follow from the results of Faudree and Schelp \cite{faudree1975path} and Gy\'{a}rf\'{a}s and Lehel \cite{gyarfas1973ramsey}. The bipartite Ramsey numbers $BR(K_{1,n},P_m)$ are given by Hatting and Henning in \cite{hattingh1998star}. The multicolor bipartite Ramsey numbers $BR(G_1,G_2,\ldots,G_t)$ when $G_1, G_2, \ldots, G_t$ are either stars and stripes, or stars and  paths, has been studied in \cite{raeisi2015star}. In \cite{bucic2019multicolour}, authors have determined asymptotically the 4-colour bipartite Ramsey number of paths and cycles. The same authors have determined asymptotically the 3-colour bipartite Ramsey number of paths and cycles in \cite{bucic20193}. The three-colour bipartite Ramsey number $BR(G_1, G_2, P_3)$ is considered in \cite{lakshmi2020three}. New values of the bipartite Ramsey number $BR(C_4,K_{1,n})$ using induced subgraphs of the incidence graph of a projective plane are given in \cite{hatala2021new}. Bipartite Ramsey numbers of  $K_{t,s}$ in many colors and bipartite Ramsey numbers of cycles for random graphs have been discussed in \cite{wang2021bipartite} and \cite{liu2021bipartite} respectively. Xuemei Zhang et.al have done a research on multicolor bipartite Ramsey numbers for quadrilaterals and stars \cite{zhang2023multicolor}. The exact value of $BR(nP_2, nP_2,....,nP_2)$ has been obtained in \cite{qiao2021note}. The exact value of $BR(P_3,P_3,P_n)$, and $BR(P_3,P_3,...,P_n)$ for $n\geq 4r+2$, and $BR(P_4,P_4,P_n)$ for $n\geq 4$ have been computed in \cite{wang2021note}.
We intend to compute the exact values of the multicolor bipartite Ramsey numbers $BR(C_8,C_{2n})$. Actually, we prove the following theorem: 

\begin{theorem}\label{th1}
	For any $n\geq 2$, we have:
	\[BR(C_8 ,C_{2n})=\begin{cases}
	8 & n=4,\\
	n+3 &\emph{otherwise.}
	\end{cases}\]
\end{theorem}
In this paper, we are only concerned with undirected, simple, and finite graphs. We follow \cite{bondy1976graph} for terminology and notations not defined here. Let $G$ be a graph with vertex set $V(G)$ and edge set $E(G)$. The degree of a vertex $v\in V(G)$ is denoted by $\deg_G(v)$, or simply by $\deg(v)$. The neighborhood $N_G(v)$ of a vertex $v$ is the set of all vertices of $G$ adjacent to $v$ and satisfies $|N_G(v)|=\deg_G(v)$. The minimum and maximum degrees of vertices of $G$ are denoted by $\delta(G)$ and $\Delta(G)$, respectively. Let $C$ be a set of colors $\{c_1,c_2,...,c_m\}$ and $E(G)$ be the edges of a graph $G$, an edge coloring
$f:E \rightarrow C$ assigns each edge in $E(G)$ to a color in $C$. If an edge
coloring uses $k$ colors on a graph, it is known as a $k$-colored graph. 
As usual, $C_n$ stands for a cycle on $n$ vertices.  The complete bipartite graph with bipartition $(X,Y)$, where $|X|=m$ and $|Y|=n$ is denoted by $K_{m,n}$. We use $[X,Y]$ to indicate the set of edges between a bipartition $(X,Y)$ of $G$. Let $W\subseteq V(G)$ be any subset of vertices of $G$, the induced subgraph $G\langle W\rangle$ is the graph whose vertex set is $W$ and whose edge set consists of all of the edges in $E(G)$ that have both endpoints in $W$.
The complement of a graph $G$, denoted by $\overline{G}$, is the graph with the same vertices as $G$ and contains those edges which are not in $G$.
$G$ is $n$-colorable to $(G_1, G_2,\ldots, G_n)$ if there exists an $n$-edge decomposition $(H_1, H_2,\ldots, H_n)$ of $G$, where $G_i\nsubseteq H_i$ for each $i=1,2, \ldots,n.$

\section{Some basic results}
To prove our main results, namely Theorem \ref{th1}, we need to establish some preliminary results. We start with the following simple but helpful lemma:
\begin{lemma}\label{le1}
	Let $G$ be a subgraph of $K_{t,t}$ with a cycle $C:x_1y_1x_2\ldots x_ky_kx_1$ of length $2k$, where $t\geq k+1$. If $x$ and $y$ be two vertices of $G$ not in $C$, where $x_i,x_{i+1}\in N_G(y)$ and $y_i,y_{i+1}\in N_G(x)$, or $xy\in E(G)$, in which $x_i\in N_G(y)$ and $y_i\in N_G(x)$ for some $i, 1\leq i\leq k$, then $G$ has a cycle of length $2k+2$.
\end{lemma}
\begin{proof}
	Consider $C'=x_1y_1x_2\ldots y_{i-1}x_iyx_{i+1}y_ixy_{i+1}x_{i+2}\ldots x_ky_kx_1$, where $x_i,x_{i+1}\in N_G(y)$ and $y_i,y_{i+1}\in N_G(x)$. Also, consider $C''=x_1y_1x_2\ldots x_iyxy_ix_{i+1}\ldots x_ky_kx_1$, where $xy\in E(G)$, $x_i\in N_G(y)$, and $y_i\in N_G(x)$.
\end{proof}
\begin{lemma}\label{le2}
	Let $G$ be a spanning subgraph of $K_{t,t}$ with a cycle $C$ of length $2k$, where $t\geq k+1$ and $k\geq 4$. Let $x$ and $y$ be the vertices of $G$ not in $C$. Assume that $x,y$ are adjacent to at least $k-1$ vertices of $C$ where $xy\notin E(G)$, or $x$ and $y$ are adjacent to at least $\lceil \frac{k}{2}\rceil+1$ vertices of $C$ where $xy\in E(G)$. Then $G$ has a cycle of length $2k+2$. 
\end{lemma}
\begin{proof}
	If $xy\notin E(G)$, the lemma is proven by Zhang et al. in \cite{zhang2013bipartite}, hence we may assume that $x,y$ are adjacent to at least $\lceil \frac{k}{2}\rceil+1$ vertices of $C$ where $xy\in E(G)$. Without loss of generality (W.l.g.), assume that $C_{2k}=x_1y_1x_2y_2\ldots x_ky_kx_1$. In this case, it is easy to check that there is at least one $i$, $1 \leq i\leq k$, such that $xy_i$ and $x_iy$ belong to $E(G)$. Therefore by Lemma \ref{le1}, we have $C_{2k+2}\subseteq G$ and the proof is complete.
\end{proof}

In the following two theorems, the authors in \cite{rui2011bipartite} and \cite{zhang2013bipartite} have determined the exact value of the bipartite Ramsey number of $BR(C_{2m}, C_{2n})$ for $m=2,3$, respectively.
\begin{theorem}\label{the1}
	For any $n\geq 2$, we have:
	\[BR(C_4 ,C_{2n})=\begin{cases}
	5 & n=2,3,\\
	n+1 &\emph{otherwise}.
	\end{cases}\]
\end{theorem}
\begin{theorem}\label{the2}
	For any $n\geq 3$, we have:
	\[BR(C_6 ,C_{2n})=\begin{cases}
	6 & n=3,\\
	n+2 &\emph{otherwise}.
	\end{cases}\]
\end{theorem}

\begin{proposition}\label{p1}
	Let $G$ be a subgraph of $K_{8,8}$ where $K_{3,4}\subseteq G$, then either $C_8\subseteq G$ or $C_8\subseteq {\overline G}$.
\end{proposition}
\begin{proof}
	Let $(X=\{x_1,\ldots,x_8\}, Y=\{y_1,\ldots ,y_8\})$ be a bipartition of $K_{8,8}$ and $K_{3,4}\subseteq G[X_1,Y_1]$, where $X_1 =\{x_1, x_2,x_3\}$, $Y_1 =\{y_1,\ldots ,y_4\}$. Assume that $C_8\nsubseteq {\overline G}$. Thus we have $|N_G(x)\cap Y_1|\leq 1$ for each $x\in X\setminus X_1$. Otherwise, $C_8\subseteq G$. Consider $X_2=X\setminus X_1$ and $Y_2= Y_1\cup \{y_5\}$. Since $BR(C_8, C_4)=5$ and $C_8\nsubseteq {\overline G}$, we have $C_4\subseteq G[X_2,Y_2]$. Hence $y_5\in V(C_4)$; otherwise, $C_8\subseteq G$. W.l.g., we may assume that $C_4=x_4y_4x_5y_5x_4$. Therefore we have $[\{x_4,x_5\},\{y_1,y_2,y_3\}], [\{x_1,x_2,x_3\},\{y_5\}]\subseteq {\overline G} $, otherwise, $C_8\subseteq G$. If $|N_{{\overline G}}(y_5)\cap \{x_6,x_7,x_8\}|\geq 2$, then $C_8\subseteq {\overline G}[X_2, \{y_1,y_2,y_3,y_5\}]$, a contradiction. So $|N_{G}(y_5)\cap \{x_6,x_7,x_8\}|\geq 2$. Let $x_6y_5, x_7y_5\in E(G)$, hence $K_{4,3}\cong [\{x_4,x_5,x_6,x_7\},\{y_1,y_2,y_3\}]\subseteq {\overline G}$. Therefore, $|N_{{\overline G}}(y)\cap \{x_4,x_5,x_6,x_7\}|\leq 1$ for each $y\in \{y_6,y_7,y_8\}$, that is $|N_{G}(y)\cap \{x_4,x_5,x_6,x_7\}|\geq 3$ for each $y\in \{y_6,y_7,y_8\}$. So, we have $K_{3,4}\cong[\{x_1,x_2,x_3\},\{y_5,y_6,y_7,y_8\}]\subseteq {\overline G}$. Therefore for each $x\in \{x_4,x_5,x_6,x_7\}$, we have $|N_{\overline{G}}(x)\cap \{y_6,y_7,y_8\}|\leq 1$; otherwise, $C_8\subseteq {\overline G}[\{x_1,x_2,x_3,x\},\{y_5,y_6,y_7,y_8\}]$ which is a contradiction. Afterwards one can check that $C_8\subseteq G[\{x_4,x_5,x_6,x_7\},\{y_5,y_6,y_7,y_8\}]$ and the proof is complete. 
\end{proof}

\section{Proof of the main results}
In this section, we compute the exact value of the bipartite Ramsey numbers $BR(C_{8}, C_{2n})$ for $n\geq 2$. Also, we guess that for $m, n\geq 3$, we have
$BR(C_{2n},C_{2m})=\begin{cases}
m+n \quad &n=m,\\
m+n-1 \quad& n\neq m.
\end{cases}$
. It should be noted that, after the authors posted the current article on arXiv, the correctness of the conjecture was proved by Yan and Peng \cite{yan2021bipartite}.\\
In order to simplify the comprehension, let us split the proof of Theorem \ref{th1} into small parts. We begin with a simple but very useful general lower bound in the following theorem:

\begin{theorem}\cite{rui2011bipartite}\label{th4}
	We have $BR(C_{2m}, C_{2n})\geq m+n-1$ for all $m, n\geq 2$.
\end{theorem}
\begin{proof}
	Let $G_1$ and $G_2$ denote vertex-disjoint induced subgraphs of $K_{m+n-2,m+n-2}$ ($m,n\geq2$), which are isomorphic to complete bipartite graphs $K_{m-1,m+n-2}$ and $K_{n-1,m+n-2}$, respectively. Clearly, $E(K_{m+n-2,m+n-2})=E(G_1)\cup E(G_2)$. As $ C_{2m}\nsubseteq G_1$ and $C_{2n}\nsubseteq G_2$, it follows that $BR(C_{2m}, C_{2n})\geq m+n-1$, as required.
\end{proof}
\begin{lemma}\cite{rui2011bipartite}\label{le3}
	$BR(C_8, C_{4})=5$.
\end{lemma}
\begin{lemma}\cite{zhang2013bipartite}\label{le4}
	$BR(C_8, C_{6})=6$.
\end{lemma}
In the following two theorems, we determine the exact values of the bipartite Ramsey number of $BR(C_8, C_{2n})$ for $n=4,5$.
\begin{theorem}\label{th5}
	$BR(C_8, C_{8})=8$.
\end{theorem}
\begin{proof}
	To prove the lower bound, consider a bipartition $(X,Y)$ of $K_{7,7}$, where $X=\{x_1,x_2,\ldots ,x_7\}$ and $Y=\{y_1,y_2,\ldots,y_7\}$. Let $(G,\overline{G})$ be a $2$-edge-coloring of $K_{7,7}$, where $G$ is given in Figure \ref{fi1}. Hence it is easy to see that $C_8\nsubseteq G$ and $\overline{G}\cong G\setminus x_4y_4$, that is we have $C_8\nsubseteq \overline{G}$. Therefore, $BR(C_8, C_{8})\geq 8$.\\
	
	\begin{figure}[ht]
		\begin{center}
		\begin{tabular}{ccc}
			\begin{tikzpicture}
			\node [draw, circle, fill=white, inner sep=2pt, label=below:$y_1$] (y1) at (0,0) {};
			\node [draw, circle, fill=white, inner sep=2pt, label=below:$y_2$] (y2) at (1,0) {};
			\node [draw, circle, fill=white, inner sep=2pt, label=below:$y_3$] (y3) at (2,0) {};
			\node [draw, circle, fill=white, inner sep=2pt, label=below:$y_4$] (y4) at (3,0) {};
			\node [draw, circle, fill=white, inner sep=2pt, label=below:$y_5$] (y5) at (4,0) {};
			\node [draw, circle, fill=white, inner sep=2pt, label=below:$y_6$] (y6) at (5,0) {};
			\node [draw, circle, fill=white, inner sep=2pt, label=below:$y_7$] (y7) at (6,0) {};
			\
			\node [draw, circle, fill=white, inner sep=2pt, label=above:$x_1$] (x1) at (0,1) {};
			\node [draw, circle, fill=white, inner sep=2pt, label=above:$x_2$] (x2) at (1,1) {};
			\node [draw, circle, fill=white, inner sep=2pt, label=above:$x_3$] (x3) at (2,1) {};
			\node [draw, circle, fill=white, inner sep=2pt, label=above:$x_4$] (x4) at (3,1) {};
			\node [draw, circle, fill=white, inner sep=2pt, label=above:$x_5$] (x5) at (4,1) {};
			\node [draw, circle, fill=white, inner sep=2pt, label=above:$x_6$] (x6) at (5,1) {};
			\node [draw, circle, fill=white, inner sep=2pt, label=above:$x_7$] (x7) at (6,1) {};
			\draw (x1)--(y1)--(x1)--(y2)--(x1)--(y3);
			\draw (x2)--(y1)--(x2)--(y2)--(x2)--(y3);
			\draw (x3)--(y1)--(x3)--(y2)--(x3)--(y3);
			\draw (x4)--(y1)--(x4)--(y2)--(x4)--(y3)--(x4)--(y4)--(x4)--(y5)--(x4)--(y6)--(x4)--(y7);
			\draw (x5)--(y5)--(x5)--(y6)--(x5)--(y7);
			\draw (x6)--(y5)--(x6)--(y6)--(x6)--(y7);
			\draw (x7)--(y5)--(x7)--(y6)--(x7)--(y7);
			\end{tikzpicture}\\
			$G$
		\end{tabular}\\
		\caption{Edge disjoint subgraphs $G$ and $\overline{G}$ of $K_{7,7}$}
		\label{fi1}
	\end{center}
	\end{figure}
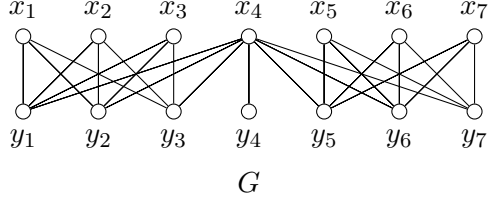
	\bigskip
	To complete the proof, suppose on the contrary that $BR(C_8, C_{8})>8$, that is $K_{8,8}$ is $2$-colorable to $(C_{8}, C_8)$, say $C_8\nsubseteq G$ and $C_{8}\nsubseteq{\overline G}$ for some $G\subseteq K_{8,8}$. Since $BR(C_6,C_{8})=6$ and $C_{8}\nsubseteq{\overline G}$, $G$ has a subgraph $C\cong C_6$. Let $(X,Y)$ be a bipartition of $K_{8,8}$, where $X=\{x_1,x_2,\ldots ,x_8\}$ and $Y=\{y_1,y_2,\ldots,y_8\}$. Set $X_1 =V(C)\cap X$, $Y_1 =V(C)\cap Y$. W.l.g., we may assume that $X_1 =\{x_1,x_2,x_3\}$, $Y_1 =\{y_1,y_2,y_3\}$, and $C_6 =x_1y_1x_2y_2x_3y_3x_1$. Consider $X_2 =X\setminus \{x_1,x_2\}$ and $Y_2=Y\setminus \{y_1,y_2\}$. Since $BR(C_6, C_{8})=6$, $|X_2|=|Y_2|=6$, and $C_{8}\nsubseteq{\overline G}$, $G[X_2,Y_2]$ has a subgraph $C'\cong C_6$. Let $X'_1 =V(C')\cap X_2$ and $Y'_1 =V(C')\cap Y_2$. Now we consider the $|\{x_3,y_3\}\cap V(C')|$, there are three cases as follows:
	
	\bigskip
	Case 1. $|\{x_3,y_3\}\cap V(C')|=2$.\\
	W.l.g., we may assume that $X'_1 =\{x_3,x_4,x_5\}$ and $Y'_1 =\{y_3,y_4,y_5\}$. If $x_3y_3\notin E(C')$, then we have $\{x_3y_4,x_3y_5, y_3x_4,y_3x_5\}\subseteq E(C')$ and there is at least one edge $x'y'$ between $\{x_4,x_5\}$ and $\{y_4,y_5\}$ in $G$. Since $x_3y_3\in E(C)$ and $x_3y', y_3x'$ are belong to $E(C')$, by Lemma \ref{le1} we have $C_8\subseteq G[X_1\cup\{x'\}, Y_1\cup\{y'\}]$, which is a contradiction. Hence, assume that $x_3y_3\in E(C')$, and w.l.g. assume that $C' =x_3y_3x_4y_4x_5y_5x_3$. Now since $C_8\nsubseteq G$, and by Figure \ref{fi2}, it is easy to check that \[\{x_1y_2,x_1y_4,x_2y_4,x_2y_5,x_4y_1,x_4y_5,x_5y_1,x_5y_2\}\subseteq E(\overline{G})\]
	That is, $C_8 =x_1y_2x_5y_1x_4y_5x_2y_4x_1\subseteq \overline{G}$, a contradiction again.
	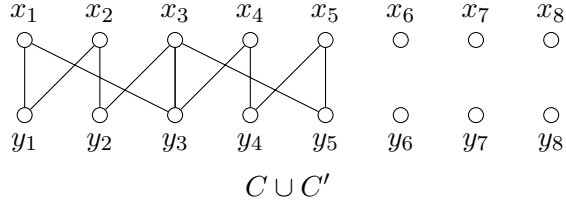
\begin{figure}[ht]
		\begin{center}
		\begin{tabular}{ccc}
			\begin{tikzpicture}
			\node [draw, circle, fill=white, inner sep=2pt, label=below:$y_1$] (y1) at (0,0) {};
			\node [draw, circle, fill=white, inner sep=2pt, label=below:$y_2$] (y2) at (1,0) {};
			\node [draw, circle, fill=white, inner sep=2pt, label=below:$y_3$] (y3) at (2,0) {};
			\node [draw, circle, fill=white, inner sep=2pt, label=below:$y_4$] (y4) at (3,0) {};
			\node [draw, circle, fill=white, inner sep=2pt, label=below:$y_5$] (y5) at (4,0) {};
			\node [draw, circle, fill=white, inner sep=2pt, label=below:$y_6$] (y6) at (5,0) {};
			\node [draw, circle, fill=white, inner sep=2pt, label=below:$y_7$] (y7) at (6,0) {};
			\node [draw, circle, fill=white, inner sep=2pt, label=below:$y_8$] (y8) at (7,0) {};
			\
			\node [draw, circle, fill=white, inner sep=2pt, label=above:$x_1$] (x1) at (0,1) {};
			\node [draw, circle, fill=white, inner sep=2pt, label=above:$x_2$] (x2) at (1,1) {};
			\node [draw, circle, fill=white, inner sep=2pt, label=above:$x_3$] (x3) at (2,1) {};
			\node [draw, circle, fill=white, inner sep=2pt, label=above:$x_4$] (x4) at (3,1) {};
			\node [draw, circle, fill=white, inner sep=2pt, label=above:$x_5$] (x5) at (4,1) {};
			\node [draw, circle, fill=white, inner sep=2pt, label=above:$x_6$] (x6) at (5,1) {};
			\node [draw, circle, fill=white, inner sep=2pt, label=above:$x_7$] (x7) at (6,1) {};
			\node [draw, circle, fill=white, inner sep=2pt, label=above:$x_8$] (x8) at (7,1) {};
			\draw (x1)--(y1)--(x2)--(y2)--(x3)--(y3)--(x1);
			\draw (x3)--(y3)--(x4)--(y4)--(x5)--(y5)--(x3);
			\end{tikzpicture}\\
			$C\cup C'$
		\end{tabular}\\
		\caption{$ x_3y_3\in E(C')$}
		\label{fi2}
		\end{center}
	\end{figure}
	
	\bigskip
	Case 2. $|\{x_3,y_3\}\cap V(C')|=1$.\\
	W.l.g., we may assume that $x_3\in V(C')$, $X'_1 =\{x_3,x_4,x_5\}$, $Y'_1 =\{y_4,y_5,y_6\}$, and $C'=x_3y_4x_4y_5x_5y_6x_3$. Since $C_8\nsubseteq G$, by Figure \ref{fi3}, it is easy to check that $K_{2,3}\cong [\{x_4,x_5\},\{y_1,y_2,y_3\}] \subseteq \overline{G}$, $K_{2,3}\cong [\{x_1,x_2\},\{y_4,y_5,y_6\}] \subseteq \overline{G}$. Now we have the following claim:
	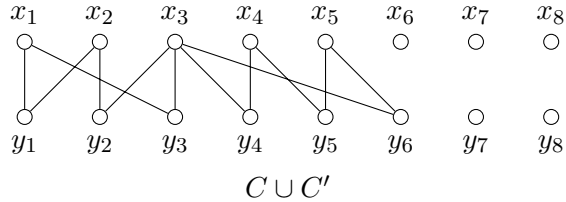
\begin{figure}[ht]
		\begin{center}
		\begin{tabular}{ccc}
			\begin{tikzpicture}
			\node [draw, circle, fill=white, inner sep=2pt, label=below:$y_1$] (y1) at (0,0) {};
			\node [draw, circle, fill=white, inner sep=2pt, label=below:$y_2$] (y2) at (1,0) {};
			\node [draw, circle, fill=white, inner sep=2pt, label=below:$y_3$] (y3) at (2,0) {};
			\node [draw, circle, fill=white, inner sep=2pt, label=below:$y_4$] (y4) at (3,0) {};
			\node [draw, circle, fill=white, inner sep=2pt, label=below:$y_5$] (y5) at (4,0) {};
			\node [draw, circle, fill=white, inner sep=2pt, label=below:$y_6$] (y6) at (5,0) {};
			\node [draw, circle, fill=white, inner sep=2pt, label=below:$y_7$] (y7) at (6,0) {};
			\node [draw, circle, fill=white, inner sep=2pt, label=below:$y_8$] (y8) at (7,0) {};
			\
			\node [draw, circle, fill=white, inner sep=2pt, label=above:$x_1$] (x1) at (0,1) {};
			\node [draw, circle, fill=white, inner sep=2pt, label=above:$x_2$] (x2) at (1,1) {};
			\node [draw, circle, fill=white, inner sep=2pt, label=above:$x_3$] (x3) at (2,1) {};
			\node [draw, circle, fill=white, inner sep=2pt, label=above:$x_4$] (x4) at (3,1) {};
			\node [draw, circle, fill=white, inner sep=2pt, label=above:$x_5$] (x5) at (4,1) {};
			\node [draw, circle, fill=white, inner sep=2pt, label=above:$x_6$] (x6) at (5,1) {};
			\node [draw, circle, fill=white, inner sep=2pt, label=above:$x_7$] (x7) at (6,1) {};
			\node [draw, circle, fill=white, inner sep=2pt, label=above:$x_8$] (x8) at (7,1) {};
			\draw (x1)--(y1)--(x2)--(y2)--(x3)--(y3)--(x1);
			\draw (x3)--(y4)--(x4)--(y5)--(x5)--(y6)--(x3);
			\end{tikzpicture}\\
			$C\cup C'$
		\end{tabular}\\
		\caption{$|\{x_3,y_3\}\cap V(C')|=1$ $(x_3\in V(C'))$}
		\label{fi3}
		\end{center}
	\end{figure}
	
	\bigskip
	\noindent\textbf{Claim 1.} If there exists a vertex $x$ of $X'_2=\{x_6,x_7,x_8\}$, such that $|N_G(x)\cap \{y_2, y_3\}|\neq 0$, then $|N_G(x)\cap Y'_1|=0$. \\
	\begin{proof} Let $xy_2\in E(G)$. For other cases, the proof is the same. Hence if $xy_4$ or $xy_6$ belong to $E(G)$, then we have $C_8\subseteq G$, and if $xy_5\in E(G)$, then we set $C' =x_3y_2xy_5x_4y_4x_3$ and the proof is the same as case $1$ when $Y_2=Y\setminus \{y_1,y_3\}$.
		\end{proof}\\
		Similar to the proof of Claim $1$, we have the following claim, which is easily verifiable.
	
	\bigskip
	\noindent\textbf{Claim 2.} If there exist a vertex $x$ of $X'_2=\{x_6,x_7,x_8\}$, such that $|N_G(x)\cap \{y_4, y_6\}|\neq 0$, then $|N_G(x)\cap Y_1|=0$.\\
	Now by Claim $1$, there are at least two vertices $x_6$ and $x_7$ of $X'_2=\{x_6,x_7,x_8\}$, such that $[\{x_6,x_7\},\{y_2, y_3\}]\subseteq {\overline G}$; otherwise, we have $K_{4,3}\subseteq \overline{G}[\{x_1,x_2,x_6,x_7,x_8\}, Y'_1]$ and the proof is complete by Proposition \ref{p1}. Therefore by Claim 2, we have $|N_G(x)\cap \{y_4, y_6\}|= 0$ for at least one vertex of $\{x_6,x_7\}$. Otherwise, we have $K_{4,3}\cong [\{x_4,x_5,x_6,x_7\}, Y_1]\subseteq {\overline G}$ and the proof is complete by Proposition \ref{p1}. Hence we may assume that $|N_G(x_7)\cap \{y_4, y_6\}|= 0$, that is $x_7y_4, x_7y_6\in E({\overline G})$. So it is easy to check that $x_1y_2, x_2y_3, x_4y_6, x_5y_4, x_6y_4,x_6y_5,x_6y_6 \in E(G)$; otherwise, $C_8\subseteq {\overline G}$. Hence $x_6y_1\in E(\overline{G})$; otherwise, $C_8\subseteq G$. Now $x_7y_1 \in E( G)$, if not we have $K_{4,3}\cong [\{x_4,x_5,x_6,x_7\}, Y_1]\subseteq {\overline G}$ and the proof is complete by Proposition \ref{p1}. We should also have  $x_7y_5\in E(\overline{G})$, if not $C_8\subseteq G$, that is we have $K_{3,3}\cong [\{x_4,x_5,x_6\}, Y_1]\subseteq {\overline G}$ and $K_{3,3}\cong [\{x_1,x_2,x_7\}, Y'_1]\subseteq {\overline G}$. Now consider $N(x_8)$, by Claims $1$ and $2$,   we have $|N_G(x_8)\cap \{y_2, y_3\}|=|N_G(x_8)\cap \{y_4, y_6\}|= 0$, otherwise we have $K_{4,3}\subseteq {\overline G}$ and the proof is complete by Proposition \ref{p1}. Therefore $\{y_2, y_3, y_4, y_6\}\subseteq N_{\overline{G}}(x_8)$ and one can check that $C_8=x_8y_6x_1y_5x_2y_4x_7y_2x_8\subseteq {\overline G}$ and the proof is complete.
	
	\bigskip
	Case 3. $|\{x_3,y_3\}\cap V(C')|=0$.\\
	We may assume that $X'_1 =\{x_4,x_5,x_6\}$, $Y'_1 =\{y_4,y_5,y_6\}$, and $C' =x_4y_4x_5y_5x_6y_6x_4$. In this case, for each $x\in X_1$, $x'\in X'_1$, $y\in Y_1$, and $y'\in Y'_1$, we have $|N_G(x)\cap Y'_1|\leq 1$, $|N_G(x')\cap Y_1|\leq 1$, $|N_G(y')\cap X_1|\leq 1$, and $|N_G(y)\cap X'_1|\leq 1$. Otherwise, the proof is the same as case 2. Therefore it is easy to show that the following claim is true.\\
	
	\noindent\textbf{Claim 3.} $K_{3,3}\setminus e \subseteq\overline{G}[X_1,Y'_1]$ and $K_{3,3}\setminus e \subseteq\overline{G}[X'_1,Y_1]$, in other words, there is at most one edge between $[X_1,Y'_1]$ and $[X'_1,Y_1]$ in $G$.\\
	Now by Claim 3, one can check that $[X_1,Y_1]\subseteq G$ or $[X'_1,Y'_1]\subseteq G$; otherwise, $C_8\subseteq \overline{G}$. W.l.g., we may assume that $[X_1,Y_1]\subseteq G$. Hence we have $[X_1,Y'_1]\subseteq \overline{G}$ or $[X'_1,Y_1]\subseteq \overline{G}$, if not, we have $C_8\subseteq G$, which is a contradiction. W.l.g., assume that $[X_1,Y'_1]\subseteq \overline{G}$. Now consider $x_7y_7$ and w.l.g. assume that $x_7y_7 \in E(G)$. For other cases, the proof is the same. If $|N_G(x_7)\cap Y_1|\neq 0$, by Lemma \ref{le1}, we have $|N_G(y_7)\cap X_1|=0$, hence $K_{3,4}\cong [X_1,Y'_1\cup \{y_7\}]\subseteq \overline{G}$ and the proof is complete by Proposition \ref{p1}. So let $|N_G(x_7)\cap Y_1|= 0$. Therefore by Claim 3, there exists one edge between $[X'_1,Y_1]$ in $G$; otherwise, $K_{3,4}\subseteq \overline{G}$ and the proof is the same. Assume that $x_4y_1\in E(G)$, if $|N_G(x_7)\cap Y'_1|= 0$ or $|N_G(y_7)\cap X_1|= 0$, then $K_{4,3}\cong [X_1\cup \{x_7\},Y'_1]\subseteq \overline{G}$ or $K_{3,4}\cong [X_1,Y'_1\cup\{y_7\}]\subseteq \overline{G}$, respectively and the proof is complete by Proposition \ref{p1}. Hence w.l.g. assume that $x_7y_4, x_1y_7\in E(G)$. Therefore one can check that $C_8 =x_1y_7x_7y_4x_4y_1x_2y_2x_1\subseteq G$, a contradiction again.\\
	Now by cases 1, 2, and 3, the proof is complete and the theorem holds.
\end{proof}
\begin{theorem}\label{th6}
	$BR(C_{10}, C_{8})=8$.
\end{theorem}
\begin{proof}
	The lower bound holds by Theorem \ref{th4}. To complete the proof, by contrary suppose that $BR(C_{10}, C_{8})>8$, that is $K_{8,8}$ is $2$-colorable to $(C_{10}, C_{8})$, say $C_{10}\nsubseteq G$ and $C_{8}\nsubseteq{\overline G}$ for some $G\subseteq K_{8,8}$. Since $BR(C_8, C_{8})=8$ and $C_{8}\nsubseteq{\overline G}$, $G$ has a subgraph $C\cong C_8$. Let $(X,Y)$ be a bipartition of $K_{8,8}$, where $X=\{x_1,x_2,\ldots ,x_8\}$ and $Y=\{y_1,y_2,\ldots,y_8\}$. Set $X_1 =V(C)\cap X$ and $Y_1 =V(C)\cap Y$. W.l.g., we may assume that $X_1 =\{x_1,x_2,x_3,x_4\}$, $Y_1 =\{y_1,y_2,y_3,y_4\}$, and $C =x_1y_1x_2y_2x_3y_3x_4y_4x_1$. Since $C_{8}\nsubseteq{\overline G}$, there is at least one edge other than $E(C)$ between $X_1$ and $Y_1$ in $G$, w.l.g. assume that $x_1y_2\in E(G)$ (for other cases, the proof is the same). Consider $X_2 =X\setminus \{x_1,x_2\}$ and $Y_2=Y\setminus \{y_1,y_2\}$. Since $BR(C_6, C_{8})=6$, $|X_2|=|Y_2|=6$, and $C_{8}\nsubseteq{\overline G}$, $G[X_2,Y_2]$ has a subgraph $C'\cong C_6$. Let $X'_1 =V(C')\cap X_1$, $Y'_1 =V(C')\cap Y_1$. Now consider the $(|X_1'|,|Y_1'|)$. By symmetry we note that $(|X_1'|,|Y_1'|)=(|Y_1'|,|X_1'|)$, so $(|X_1'|,|Y_1'|) \in \{(0,0),(1,0), (2,0), (1,1), (2,1), (2,2)\}$, now we have the following cases: 
	
	\bigskip
	Case 1. $(|X_1'|,|Y_1'|)=(0,0)$.\\
	Assume that $V(C') =\{x_5,x_6,x_7, y_5,y_6,y_7\}$ and $C' =x_5y_5x_6y_6x_7y_7x_5$. 
	Set $X'_2 =\{x_5,x_6,x_7\}$ and $Y'_2 =\{y_5,y_6,y_7\}$. Since $C_{10}\nsubseteq G$, one can check that either $K_{3,4}\cong [X'_2,Y_1]\subseteq \overline{G}$ or $K_{4,3} \cong [X_1,Y'_2]\subseteq \overline{G}$; otherwise, assume to the contrary that there exists at least one edge between  $X'_2$ and $Y_1$ in $G$, along with there exists at least one edge between $X_1$ and $Y'_2$ in $G$.  Utilizing symmetry and without loss of generality, let's consider that $x_5y_1\in E(G)$. Subsequently, let's posit that  $|N_G(y_5)\cap X_1|\neq 0$. For other cases, the proof remains consistent.  Since both $x_5y_5$ and $x_5y_1$ are edges in $G$, if $x_1y_5$ or $x_2y_5$ are also edges in $G$, then by lemma \ref{le1} it can be deduced tha $C_{10}\subseteq G$. Therefore, let's assume that either $x_3y_5\in E(G)$ or $x_4y_5\in E(G)$. In the instance of $x_3y_5\in E(G)$, we have $C =x_5y_1x_2y_2x_1y_4x_4y_3x_3y_5x_5$ is a copy of $C_{10}$ in $G$.  The proof remains analogous in the scenario where  $x_4y_5\in E(G)$.\\
	Hence,  w.l.g., we may assume that $K_{3,4}\cong [X'_2,Y_1]\subseteq \overline{G}$. If there exists a vertex $x$ of $X\setminus X'_2$, such that $|N_{\overline{G}}(x)\cap Y_1|\geq 2$, then $C_{8}\subseteq{\overline G}$, a contradiction. Hence $|N_G(x)\cap Y_1|\geq 3$ for each $x\in X\setminus X'_2$, that is $|N_G(x_8)\cap Y_1|\geq 3$. Therefore by Lemma \ref{le1}, we have $|N_G(y)\cap X_1|\leq 1$ for each $y\in Y\setminus Y_1 $; otherwise, $C_{10}\subseteq G$. If there exist $y,y'\in \{y_5,y_6,y_7\}$, such that $|N_G(y)\cap X_1|=|N_G(y')\cap X_1|=1$ and $N_G(y)\cap X_1\neq N_G(y')\cap X_1$,  one can check that $C_{10}\subseteq G$, a contradiction again (for example w.l.g., assume that $N_G(y_5)\cap X_1=\{x_1\}$ and  $N_G(y_6)\cap X_1=\{x_2\}$. Hence $C_{10}:=y_5x_1y_4x_4y_3x_3y_2x_2y_6x_6y_5\subseteq G$). Therefore, w.l.g. we may assume that $\{x_1,x_2,x_3\}\subseteq N_{{\overline G}}(y)$ for each $y\in\{y_5,y_6,y_7\}$. If there are at least two vertices $y_5,y_6$ of $Y'_2$, such that $|N_{{\overline G}}(y_i)\cap X_1|= 4$ for $i=5,6$, then $C_{8}\subseteq{\overline G}[X_1,Y\setminus Y_1]$. In other words, there are at least two vertices $y_5,y_6$ of $Y'_2$,  such that $N_G(y_i)\cap X_1= \{x_4\}$. Hence $x_8y_i\in E(\overline{G})$ for $i=5,6$, if not, $C_{10}\subseteq G$, a contradiction. Now one can check that $C_{8}\subseteq{\overline G}[\{x_1,x_2,x_3,x_8\},Y\setminus Y_1]$ and the proof is complete.
	
	\bigskip
	Case 2. $(|X_1'|,|Y_1'|)=(1,0)$.\\
	Assume that $X'_1=\{x_4\}$. For the case that $X'_1=\{x_3\}$, the proof is the same. Hence w.l.g. assume that $V(C') =\{x_4,x_5,x_6, y_5,y_6,y_7\}$ and $C' =x_4y_5x_5y_6x_6y_7x_4$. 
	Set $X''_2 =\{x_1,x_2,x_3\}$, $X'_2 =\{x_5,x_6\}$, and $Y'_2 =\{y_5,y_6,y_7\}$. Since $C_{10}\nsubseteq G$, we have $K_{3,3}\cong [X''_2,Y'_2]\subseteq \overline{G}$ and $K_{2,4} \cong [X'_2,Y_1]\subseteq \overline{G}$. Now consider the vertices $\{x_7, x_8\}$, one can check that for at least one vertex $x$ of $\{x_7, x_8\}$, we have  $|N_G(x)\cap Y_1|\geq 2$; otherwise, by$K_{2,4} \cong [X'_2,Y_1]\subseteq \overline{G}$ one can say that  $C_{8}\subseteq{\overline G}[X\setminus X_1,Y_1]$. W.l.g., assume that $|N_G(x_7)\cap Y_1|\geq 2$ and thus  $|N_G(x_7)\cap Y'_2|=0$; otherwise, $C_{10}\subseteq G$. Therefore, we have $K_{4,3}\cong [ X''_2\cup\{x_7\},Y'_2]\subseteq \overline{G}$, that is $|N_{{\overline G}}(y_8)\cap (X''_2\cup\{x_7\})|\leq 1$, if  not, $C_{8}\subseteq{\overline G}[X''_2\cup\{x_7\},Y_2]$. Since $|N_G(x_7)\cap Y_1|\geq 2$ and $|N_G(y_8)\cap X''_2|\geq 2$, if $x_7y_8 \in E(G)$, using  Lemma \ref{le1}, it can be easily checked that $C_{10}\subseteq G$. So let $x_7y_8 \in E({\overline G})$, that is $X''_2 \subseteq N_G(y_8)$ and thus $|N_G(x_7)\cap Y_1|=2$ and $N_G(x_7)\cap Y_1=\{y_3,y_4\}$. Otherwise by Lemma \ref{le1}, we have $C_{10}\subseteq G$ and the proof is complete. Similarly $|N_G(x_8)\cap Y_1|\leq 2$; otherwise, by Lemma \ref{le2}, $C_{10}\subseteq G[X_1\cup\{x_8\},Y_1\cup\{y_8\}]$. If $|N_G(x_8)\cap Y_1| \leq 1$, then $C_{8}\subseteq \overline{G}[X\setminus X_1,Y_1]$. Hence $|N_G(x_8)\cap Y_1|=2$ and $x_8y_8 \notin E(G)$, if not, one can by Lemma \ref{le1}, check that $C_{10}\subseteq G[X_1\cup\{x_8\},Y_1\cup\{y_8\}]$. Thus we have $Y\setminus Y_1\subseteq N_{{\overline G}}(x_i)$ for $i=7,8$. That is, $C_{8}\subseteq \overline{G}[ \{x_1,x_2,x_7,x_8\} ,Y\setminus Y_1]$ and the proof is complete.
	
	\bigskip
	Case 3. $(|X_1'|,|Y_1'|)=(1,1)$.\\
	W.l.g., we may assume that $V(C') =\{x,x_5,x_6, y,y_5,y_6\}$, where $x\in \{x_3,x_4\}$ and $y\in \{y_3,y_4\}$. If $xy\notin E(C)$, then $x=x_3$ and $y=y_4$. Let $x_3y_4\in E(C')$, and w.l.g. we may assume that $C' =x_3y_4x_5y_5x_6y_6x_3$. In this case, we have $C_{10} = x_1y_1x_2y_2x_3y_6x_6y_5x_5y_4x_1\subseteq G$, a contradiction. So let $x_3y_4\notin E(C')$ and w.l.g. we may assume that $C' =x_3y_5x_5y_4x_6y_6x_3$. Consider the following figure:\\
	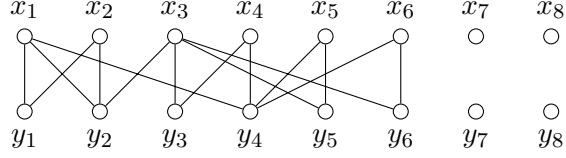
\begin{figure}[ht]
		\begin{center}
		\begin{tabular}{ccc}
			\begin{tikzpicture}
			\node [draw, circle, fill=white, inner sep=2pt, label=below:$y_1$] (y1) at (0,0) {};
			\node [draw, circle, fill=white, inner sep=2pt, label=below:$y_2$] (y2) at (1,0) {};
			\node [draw, circle, fill=white, inner sep=2pt, label=below:$y_3$] (y3) at (2,0) {};
			\node [draw, circle, fill=white, inner sep=2pt, label=below:$y_4$] (y4) at (3,0) {};
			\node [draw, circle, fill=white, inner sep=2pt, label=below:$y_5$] (y5) at (4,0) {};
			\node [draw, circle, fill=white, inner sep=2pt, label=below:$y_6$] (y6) at (5,0) {};
			\node [draw, circle, fill=white, inner sep=2pt, label=below:$y_7$] (y7) at (6,0) {};
			\node [draw, circle, fill=white, inner sep=2pt, label=below:$y_8$] (y8) at (7,0) {};
			\
			\node [draw, circle, fill=white, inner sep=2pt, label=above:$x_1$] (x1) at (0,1) {};
			\node [draw, circle, fill=white, inner sep=2pt, label=above:$x_2$] (x2) at (1,1) {};
			\node [draw, circle, fill=white, inner sep=2pt, label=above:$x_3$] (x3) at (2,1) {};
			\node [draw, circle, fill=white, inner sep=2pt, label=above:$x_4$] (x4) at (3,1) {};
			\node [draw, circle, fill=white, inner sep=2pt, label=above:$x_5$] (x5) at (4,1) {};
			\node [draw, circle, fill=white, inner sep=2pt, label=above:$x_6$] (x6) at (5,1) {};
			\node [draw, circle, fill=white, inner sep=2pt, label=above:$x_7$] (x7) at (6,1) {};
			\node [draw, circle, fill=white, inner sep=2pt, label=above:$x_8$] (x8) at (7,1) {};
			\draw (x1)--(y1)--(x2)--(y2)--(x3)--(y3)--(x4)--(y4)--(x1)--(y2);
			\draw (x3)--(y5)--(x5)--(y4)--(x6)--(y6)--(x3);
			\end{tikzpicture}\\
		\end{tabular}\\
		\caption{$x=x_3, y=y_4$ and $x_3y_4\notin E(C')$}
		\label{fi4}
		\end{center}
	\end{figure}

	By Figure \ref{fi4}, it is easy to check that $x_iy_j\in E(\overline{G})$, where $i\in\{5,6\}$, $j\in\{1,2,3\}$. Similarly, $x_jy_i\in E(\overline{G})$, where $i\in\{5,6\}$, $j\in\{1,2,4\}$, and $x_2y_3, x_4y_2 \in E(\overline{G})$. Hence we have $C_8\subseteq \overline{G}[\{x_2,x_4,x_5,x_6\},\{y_1,y_2,y_3,y_5\}]$, a contradiction too. So let $xy\in E(C)$, that is $xy\in \{x_3y_3,y_3x_4,x_4y_4\}$. In this case, the proof is the same as the case that $x=x_3$ and $y=y_4$, where $x_3y_4\in E(C')$ and we get a contradiction again.
	
	\bigskip
	Case 4. $(|X_1'|,|Y_1'|)=(2,0)$.\\
	W.l.g., let $V(C') =\{x_3,x_4,x_5, y_5,y_6,y_7\}$ and $C' =x_3y_5x_4y_6x_5y_7x_3$. In this case, we have $C_{10} = x_1y_1x_2y_2x_3y_7x_5y_6x_4y_4x_1\subseteq G$,  which is a contradiction.
	
	\bigskip
	Case 5. $(|X_1'|,|Y_1'|)=(2,1)$.\\
	W.l.g., we may assume that $V(C') =\{x_3,x_4,x_5, y,y_5,y_6\}$, where $y\in \{y_3,y_4\}$. If $y=y_3$ or $y=y_4$ and  $x_3y_4\in E(C')$,  by considering the edges of $C$ and $C'$, in any case, it is easy to check that $C_{10}\subseteq G[V(C)\cup V(C')]$. For example, assume that $y=y_3$ and $C'=x_3y_3x_4y_5x_5y_6x_3$, hence we have $C_{10}:=y_6x_5y_5x_4y_4x_1y_1x_2y_2x_3y_6\subseteq G$. For other cases, the proof is the same.
	
	Hence, assume that $y=y_4$, and $x_3y_4\notin E(C')$. W.l.g., assume that $C'=x_4y_4x_5y_5x_3y_6x_4$. Consider the following figure:\\
	\begin{figure}[ht]
		\begin{center}
		\begin{tabular}{ccc}
			\begin{tikzpicture}
			\node [draw, circle, fill=white, inner sep=2pt, label=below:$y_1$] (y1) at (0,0) {};
			\node [draw, circle, fill=white, inner sep=2pt, label=below:$y_2$] (y2) at (1,0) {};
			\node [draw, circle, fill=white, inner sep=2pt, label=below:$y_3$] (y3) at (2,0) {};
			\node [draw, circle, fill=white, inner sep=2pt, label=below:$y_4$] (y4) at (3,0) {};
			\node [draw, circle, fill=white, inner sep=2pt, label=below:$y_5$] (y5) at (4,0) {};
			\node [draw, circle, fill=white, inner sep=2pt, label=below:$y_6$] (y6) at (5,0) {};
			\node [draw, circle, fill=white, inner sep=2pt, label=below:$y_7$] (y7) at (6,0) {};
			\node [draw, circle, fill=white, inner sep=2pt, label=below:$y_8$] (y8) at (7,0) {};
			\
			\node [draw, circle, fill=white, inner sep=2pt, label=above:$x_1$] (x1) at (0,1) {};
			\node [draw, circle, fill=white, inner sep=2pt, label=above:$x_2$] (x2) at (1,1) {};
			\node [draw, circle, fill=white, inner sep=2pt, label=above:$x_3$] (x3) at (2,1) {};
			\node [draw, circle, fill=white, inner sep=2pt, label=above:$x_4$] (x4) at (3,1) {};
			\node [draw, circle, fill=white, inner sep=2pt, label=above:$x_5$] (x5) at (4,1) {};
			\node [draw, circle, fill=white, inner sep=2pt, label=above:$x_6$] (x6) at (5,1) {};
			\node [draw, circle, fill=white, inner sep=2pt, label=above:$x_7$] (x7) at (6,1) {};
			\node [draw, circle, fill=white, inner sep=2pt, label=above:$x_8$] (x8) at (7,1) {};
			\draw (x1)--(y1)--(x2)--(y2)--(x3)--(y3)--(x4)--(y4)--(x1)--(y2);
			\draw (x3)--(y5)--(x5)--(y4)--(x4)--(y6)--(x3);
			\end{tikzpicture}\\
		\end{tabular}\\
		\caption{$x_3,x_4 \in X_1',y_4\in Y_1'$ and $x_3y_4\notin E(C')$}
		\label{fi5}
		\end{center}
	\end{figure}
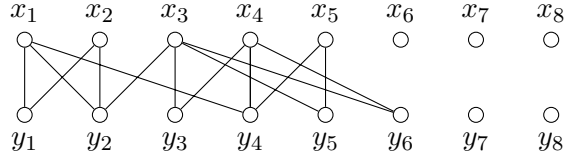

	By Figure \ref{fi5}, it is easy to check that:
	\[\{x_1y_3,x_1y_5,x_1y_6,x_2y_3,x_2y_5,x_2y_6,x_4y_5, x_5y_2, x_5y_3, x_5y_6 \}\subseteq E(\overline{G})\]
	Otherwise, $C_{10}\subseteq G$.  For example by contrary assume that $x_1y_3\in E(G)$, it can be seen that $C_{10}:=x_1y_1x_2y_2x_3y_5x_5y_4x_4y_3x_1\subseteq G$. For other cases, the proof is the same. Now, consider $x_4y_2$. If $x_4y_2\in E(\overline{G})$, then we have $C_8=x_1y_3x_5y_2x_4y_5x_2y_6x_1\subseteq \overline{G}$, a contradiction. Hence assume that $x_4y_2\in E(G)$, therefore $C_{10}=x_1y_1x_2y_2x_4y_6x_3y_5x_5y_4x_1\subseteq G$, which is a contradiction again.
	
	\bigskip
	Case 6. $(|X_1'|,|Y_1'|)=(2,2)$.\\
	W.l.g., we may assume that $V(C') =\{x_3,x_4,x_5, y_3,y_4,y_5\}$. If $x_5y_5\notin E(C')$, we have $x_5y_j , y_5x_j \in E(C')$ for $j=3,4$. Thus by Lemma \ref{le1}, we have $C_{10}\subseteq G$, which is a contradiction. Now let $x_5y_5\in E(C')$. If $x_4y_5\in E(C')$, then by Lemma \ref{le1}, the proof is the same. So let $x_4y_5\notin E(C')$. Therefore, we have $x_3y_5\in E(C')$. If $x_5y_3\in E(C')$, the proof is the same. Hence $C' =x_3y_3x_4y_4x_5y_5x_3 \subseteq G$. Since $C_{10}\nsubseteq G$, and by Figure \ref{fi6}, it is easy to check that:
	\[\{x_1y_3,x_1y_5,x_2y_3,x_2y_5, x_4y_1,x_4y_2,x_4y_5,x_5y_1,x_5y_2,x_5y_3\}\subseteq E(\overline{G}).\]\\
	\begin{figure}[ht]
		\begin{center}
		\begin{tabular}{ccc}
			\begin{tikzpicture}
			\node [draw, circle, fill=white, inner sep=2pt, label=below:$y_1$] (y1) at (0,0) {};
			\node [draw, circle, fill=white, inner sep=2pt, label=below:$y_2$] (y2) at (1,0) {};
			\node [draw, circle, fill=white, inner sep=2pt, label=below:$y_3$] (y3) at (2,0) {};
			\node [draw, circle, fill=white, inner sep=2pt, label=below:$y_4$] (y4) at (3,0) {};
			\node [draw, circle, fill=white, inner sep=2pt, label=below:$y_5$] (y5) at (4,0) {};
			\node [draw, circle, fill=white, inner sep=2pt, label=below:$y_6$] (y6) at (5,0) {};
			\node [draw, circle, fill=white, inner sep=2pt, label=below:$y_7$] (y7) at (6,0) {};
			\node [draw, circle, fill=white, inner sep=2pt, label=below:$y_8$] (y8) at (7,0) {};
			\
			\node [draw, circle, fill=white, inner sep=2pt, label=above:$x_1$] (x1) at (0,1) {};
			\node [draw, circle, fill=white, inner sep=2pt, label=above:$x_2$] (x2) at (1,1) {};
			\node [draw, circle, fill=white, inner sep=2pt, label=above:$x_3$] (x3) at (2,1) {};
			\node [draw, circle, fill=white, inner sep=2pt, label=above:$x_4$] (x4) at (3,1) {};
			\node [draw, circle, fill=white, inner sep=2pt, label=above:$x_5$] (x5) at (4,1) {};
			\node [draw, circle, fill=white, inner sep=2pt, label=above:$x_6$] (x6) at (5,1) {};
			\node [draw, circle, fill=white, inner sep=2pt, label=above:$x_7$] (x7) at (6,1) {};
			\node [draw, circle, fill=white, inner sep=2pt, label=above:$x_8$] (x8) at (7,1) {};
			\draw (x1)--(y1)--(x2)--(y2)--(x3)--(y3)--(x4)--(y4)--(x1)--(y2);
			\draw (x3)--(y3)--(x4)--(y4)--(x5)--(y5)--(x3);
			\end{tikzpicture}\\
		\end{tabular}\\
		\caption{$(|X_1'|,|Y_1'|)=(2,2)$ and $x_4y_5, x_5y_3 \in E(\overline{G})$}
		\label{fi6}
		\end{center}
	\end{figure}
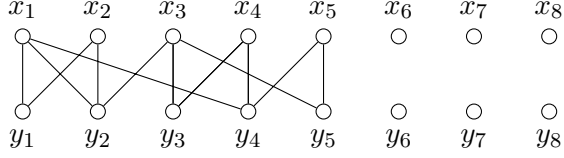

	Now we have the following claim:
	
	\bigskip
	\noindent\textbf{Claim 4.} $|N_G(x_5)\cap \{y_6,y_7,y_8\}|=0$. \\
	\begin{proof} By contradiction we may assume that $|N_G(x_5)\cap \{y_6,y_7,y_8\}|\neq 0$, say $x_5y_6\in E(G)$. In this case, we have $\{x_1,x_2,x_4\}\subseteq N_{\overline{G}}(y_6)$; otherwise, $C_{10}\subseteq G$ which is a contradiction. Now, since $\{x_1,x_2,x_4\}\subseteq N_{\overline{G}}(y_6)$, one can check that $C_8 =y_6x_4y_1x_5y_3x_1y_5x_2y_6\subseteq \overline{G}$, a contradiction again. Hence the assumption does not hold and the claim follows.\end{proof}\\
	Therefore by Claim $4$, $\{y_6,y_7,y_8\}\subseteq N_{\overline{G}}(x_5)$, and we have the following claim:
	
	\bigskip
	\noindent\textbf{Claim 5.} $\{x_1,x_2\}\subseteq N_{G}(y_i) $ for $i=6,7,8$. \\
	\begin{proof}
	 By contradiction, let there is at least one edge between $\{x_1,x_2\}$ and $\{y_6,y_7,y_8\}$ in $\overline{G}$, say $x_1y_6\in E(\overline{G})$. Hence $C_8 = x_1y_3x_2y_5x_4y_1x_5y_6x_1\subseteq \overline{G}$, a contradiction.\\
	Now by Lemma \ref{le1} and Claim $4$, for each $i=6,7,8$, $x_4y_i\in E(\overline{G})$; otherwise, $C_{10}\subseteq G$, a contradiction. If there exists a vertex $x$ of $\{x_6,x_7,x_8\}$, such that $|N_{\overline{G}}(x)\cap \{y_6,y_7,y_8\}|\geq 2$, then we have $C_8\subseteq \overline{G}$, a contradiction as well. Hence, for each $i=6,7,8$, we have $|N_G(x)\cap \{y_6,y_7,y_8\}|\geq 2$. Assume that $x_6y_6,x_6y_7\in E(G)$, therefore it is easy to check that $C_{10}\subseteq G[X_1\cup\{x_6\}, Y_1\setminus\{y_1\}\cup\{y_6,y_7\}]$, a contradiction again.
	\end{proof}
	
	Now by cases $1,2,\ldots,6$, the proof is complete and the theorem holds.
\end{proof}

In the following theorem, we determine the exact value of the bipartite Ramsey number  $BR(C_8, C_{2n})$ for $n\geq 5$.
\begin{theorem}\label{th7}
	$BR(C_8,C_{2n})=n+3$ for each $n\geq 5$.
\end{theorem}
\begin{proof}
	The lower bound holds by Theorem \ref{th4}. We use induction to prove the upper bound. For the base step of the induction, the theorem holds by Theorem \ref{th6}. Suppose that $n\geq 6$ and $BR(C_8,C_{2n'})=n'+3$ for each $n'<n$. We will show that $BR(C_8,C_{2n})\leq n+3$. To complete the proof, by contrary suppose that $BR(C_8,C_{2n})>n+3$, that is there exists a subgraph $G$ of $K_{t,t}$, such that neither $C_{2n}\subseteq G$ nor $C_{8}\subseteq{\overline G}$, where $t=n+3$. Since $BR(C_8,C_{2(n-1)})=n+2$ and $C_{8}\nsubseteq{\overline G}$, $G$ has a subgraph $C\cong C_{2(n-1)}$. Let $(X,Y)$ be a bipartition of $K_{t,t}$, where $X=\{x_1,x_2,\ldots ,x_t\}$ and $Y=\{y_1,y_2,\ldots,y_t\}$. Set $X_1 =V(C)\cap X$ and $Y_1 =V(C)\cap Y$. W.l.g., we may assume that $X_1 =\{x_1,x_2,\ldots x_{n-2},x_{n-1}\}$, $Y_1 =\{y_1,y_2,\ldots,y_{n-2},y_{n-1}\}$,  and $C =x_1y_1x_2y_2\ldots x_{n-2}y_{n-2}x_{n-1}y_{n-1}x_1$. Since $C_{8}\nsubseteq{\overline G}$, there is at least one edge other than $E(C)$ between $X_1$ and $Y_1$ in $G$. Now we have the following cases.
	
	\bigskip 
	Case 1. There exists $x_iy_k\in E(G)$ for some $i\in\{1,2,\ldots,n-1\}$, where $k-i=1(\mod n-1)$ or $i-k=2(\mod n-1)$.\\
	W.l.g., assume that $x_1y_2\in E(G)$. For other cases, the proof is the same. 
	Set $X_2 = \{x_{n-2},x_{n-1}, \ldots,x_{n+3}\}$, $Y_2=\{y_{n-2},y_{n-1}, \ldots,y_{n+3}\}$, $X''=\{x_1,x_2,x_{n-2},x_{n-1}\}$, and  $Y''=\{y_1,y_2,y_{n-2},y_{n-1}\}$. Since $BR(C_6, C_{8})=6$, $|X_2|=|Y_2|=6$ and $C_{8}\nsubseteq{\overline G}$, $G[X_2,Y_2]$ has a subgraph $C'\cong C_6$. Let $X'_1 =V(C')\cap X_1$ and $Y'_1 =V(C')\cap Y_1$. Now consider the $(|X_1'|,|Y_1'|)$, we note that $(|X_1'|,|Y_1'|)=(|Y_1'|,|X_1'|)$, hence we have $(|X_1'|,|Y_1'|)\in \{(0,0),(1,0), (2,0), (1,1), (2,1), (2,2)\}$. Now we have the following subcases:
	
	\bigskip 
	Subcase 1-1. $(|X_1'|,|Y_1'|)=(0,0)$.\\
	W.l.g. assume that $V(C') =\{x_{n},x_{n+1},x_{n+2}, y_{n},y_{n+1},y_{n+2}\}$ and $C' =x_{n}y_{n}x_{n+1}y_{n+1}x_{n+2}y_{n+2}x_{n}$. Set $X'_2 =\{x_{n},x_{n+1},x_{n+2}\}$, $Y'_2 =\{y_{n},y_{n+1},y_{n+2}\}$. Since $C_{2n}\nsubseteq G$, one can check that either $K_{3,4}\cong [X'_2,Y'']\subseteq \overline{G}$ or $K_{4,3} \cong [X'',Y'_2]\subseteq \overline{G}$. Otherwise, let's assume the contrary: there is at least one edge between the sets $X'_2$ and $Y''$ within the graph $G$, and there exists at least  one edge between $X''$ and $Y'_2$ in $G$. Without loss of generality and due to symmetry, let's assume that $x_{n}y_1\in E(G)$. Now, let's assume that $|N_G(y_n)\cap X_1|\neq 0$.  The proof remains the same for other cases. As $x_ny_n\in E(G)$ and $x_ny_1\in E(G)$, if either $x_1y_n$ or $x_2y_n\in E(G)$, then, according to Lemma \ref{le1}  it can be concluded that $C_{2n}\subseteq G$. Therefore, let's suppose that either $x_{n-2}y_5n\in E(G)$ or $x_{-1}ny_n\in E(G)$. Assuming $x_{n-2}y_n\in E(G)$,  it follows that  $C =x_ny_1x_2y_2x_3y_3\ldots x_{n-3}y_{n-3}x_{n-2}y_nx_{n+1}y_{n+1}x_{n+2}y_{n+2}x_n$ forms a copy of $C_{2n}$ in $G$. The proof remains the same in the scenario where  $x_{n-1}y_n\in E(G)$.\\
	Now w.l.g. let $K_{3,4}\cong [X'_2,Y''\}]\subseteq \overline{G}$. So for each $x\in X\setminus X'_2$, we have $|N_{\overline{G}}(x)\cap Y''|\leq 1$; otherwise, $C_{8}\subseteq{\overline G}$, which is a contradiction. That is, $|N_G(x)\cap Y''|\geq 3$ for each $x\in X''$ and $|N_G(x_{n+3})\cap Y''|\geq 3$. Thus, one can assume that $|N_G(y)\cap X''|\leq 1$ for each $y\in Y\setminus Y_1$.  If this were not the case, assume by contradiction that  $|N_G(y)\cap X''|\geq 2$ for at least one $y\in Y\setminus Y_1$.  Considering the situation, where  $|N_G(x)\cap Y''|\geq 3$ holds true for every $x\in X''$ and $|N_G(x_{n+3})\cap Y''|\geq 3$,  it becomes evident that there is an $i\in \{1,2,n-2,n-1\}$, such that $x_iy_ix_{i+1}y_{i+1}$ forms a part of a copy of  $C_{2n-2}$ within $G[X_1,Y_1]$. In this configuration, $x_i,x_{i+1}\in N_G(y)$ and $y_i, y_{i+1}\in N_G(x_{n+3})$. By utilizing Lemma  \ref{le1},  it is concluded that $C_{2n}\subseteq G$, leading to a contradiction. If there exist $y,y'\in \{y_n,y_{n+1},y_{n+2}\}$, such that $|N_G(y)\cap X''|=|N_G(y')\cap X''|=1$ and $N_G(y)\cap X''\neq N_G(y')\cap X''$, then one can check that $C_{2n}\subseteq G$, a contradiction too. W.l.g., we may assume that $\{x_1,x_2,x_{n-2}\}\subseteq N_{\overline{G}}(y)$ for each $y\in \{y_n,y_{n+1},y_{n+2}\}$. If there are at least two vertices of $\{y_n,y_{n+1},y_{n+2}\}$, say $y_n,y_{n+1}$, such that $|N_{{\overline G}}(y_i)\cap X''|= 4$ for $i=n,n+1$, then $C_{8}\subseteq{\overline G}[X'',Y\setminus Y_1]$. In other words, there are at least two vertices of $\{y_n,y_{n+1},y_{n+2}\}$, say $y_n,y_{n+1}$, such that $N_G(y_i)\cap X''= \{x_{n-1}\}$. Therefore $x_{n+3}y_i\in E(\overline{G})$ for $i=n,n+1$; otherwise, $C_{2n}\subseteq G$, a contradiction again. Now one can check that $C_{8}\subseteq{\overline G}[\{x_1,x_2,x_{n-2},x_{n+3}\},Y\setminus Y_1]$ and the proof is complete.
	
	\bigskip
	Subcase 1-2. $(|X_1'|,|Y_1'|)=(1,0)$.\\
	W.l.g., assume that $V(C') =\{x_{n-1},x_{n},x_{n+1}, y_{n},y_{n+1},y_{n+2}\}$ and $C' =x_{n-1}y_{n}x_{n}y_{n+1}x_{n+1}y_{n+2}x_{n-1}$. Set $X''_2 =\{x_1,x_2,x_{n-2}\}$, $X'_2 =\{x_n,x_{n+1}\}$, and $Y'_2 =\{y_n,y_{n+1},y_{n+2}\}$.\\ Hence, one can assume that $K_{3,3}\cong [X''_2,Y'_2]\subseteq \overline{G}$. Otherewise, let $|N_G(y)\cap X''_2|\neq 0$ for at least one $y\in Y'_2$. W.l.g, assume that  $x_1y_n\in E(G)$, the proof for other cases follows similarly. In this case, we observe that the cycle $C =y_nx_1y_2x_3y_3\ldots x_{n-2}y_{n-2}x_{n-1}y_{n+2}x_{n+1}y_{n+1}x_{n}y_n$ forms a copy of  $C_{2n}$ within  $G$, which is a contradiction. Applying symmetry, we can similarly deduce that $K_{2,4} \cong [X'_2,Y'']\subseteq \overline{G}$.\\
	Now consider the vertices $\{x_{n+2}, x_{n+3}\}$. One can check that $|N_G(x)\cap Y''|\geq 2$ for at least one $x\in \{x_{n+2}, x_{n+3}\}$; otherwise, we have $C_{8}\subseteq{\overline G}[\{x_{n},x_{n+1},x_{n+2},x_{n+3}\},Y'']$. W.l.g., we may assume that $|N_G(x_{n+2})\cap Y''|\geq 2$, hence we have $|N_G(x_{n+2})\cap Y'_2|=0$, if not, $C_{2n}\subseteq G$. Therefore $K_{4,3}\cong [ X''_2\cup\{x_{n+2}\},Y'_2]\subseteq \overline{G}$, and so $|N_{{\overline G}}(y_{n+3})\cap (X''_2\cup\{x_{n+2}\})|\leq 1$; otherwise, $C_{8}\subseteq{\overline G}[X''_2\cup\{x_{n+2}\},Y'_2\cup\{y_{n+3}\}]$. Since $|N_G(x_{n+2})\cap Y''|\geq 2$ and $|N_G(y_{n+3})\cap X''_2|\geq 2$, if $x_{n+2}y_{n+3} \in E(G)$ then by Lemma \ref{le1}, we have $C_{2n}\subseteq G$. Hence $x_{n+2}y_{n+3} \in E({\overline G})$, that is $X''_2 \subseteq N_G(y_{n+3})$. Thus $|N_G(x_{n+2})\cap Y''|=2$ and $N_G(x_{n+2})\cap Y''=\{y_{n-2},y_{n-1}\}$; otherwise, by Lemma \ref{le1}, we have $C_{2n}\subseteq G$ and the proof is complete. Similarly $|N_G(x_{n+3})\cap Y''|\leq 2$, if not, by Lemma \ref{le1}, $C_{2n}\subseteq G[X_1\cup\{x_{n+3}\},Y_1\cup\{y_{n+3}\}]$. If $|N_G(x_{n+3})\cap Y''| \leq 1$, we have $C_{8}\subseteq \overline{G}[X\setminus X_1,Y'']$. So $|N_G(x_{n+3})\cap Y''|=2$ and $x_{n+3}y_{n+3} \notin E(G)$; otherwise, by Lemma \ref{le1}, it is easy to check that $C_{2n}\subseteq G[X_1\cup\{x_{n+3}\},Y_1\cup\{y_{n+3}\}]$. Thus we have $Y\setminus Y_1\subseteq N_{{\overline G}}(x_i)$ for $i=n+2,n+3$. That is, we have $C_{8}\subseteq \overline{G}[ \{x_1,x_2,x_{n+2},x_{n+3}\} ,Y\setminus Y_1]$, and the proof is complete.
	
	\bigskip
	Subcase 1-3. $(|X_1'|,|Y_1'|)=(1,1)$.\\
	W.l.g., assume that $V(C') =\{x,x_{n},x_{n+1}, y, y_{n},y_{n+1}\}$, where $x\in \{x_{n-2},x_{n-1}\}$ and $ y\in \{y_{n-2},y_{n-1}\} $. If $xy\notin E(C)$, then $x=x_{n-2}$ and $y=y_{n-1}$. Let $x_{n-2}y_{n-1}\in E(C')$ and assume that $C' =x_{n-2}y_{n-1}x_ny_nx_{n+1}y_{n+1}x_{n-2}$. In this case, we have $C_{2n} = x_1y_1x_2y_2\ldots x_{n-2}y_{n+1}x_{n+1}y_nx_ny_{n-1}x_1\subseteq G$, a contradiction. Assume that $x_{n-2}y_{n-1}\notin E(C')$ and w.l.g. let $C' =x_{n-2}y_nx_ny_{n-1}x_{n+1}y_{n+1}x_{n-2}$. By considering Figure \ref{fi7}, it is easy to check that: 
	\[\{x_2y_{n-2}, x_2y_n, x_{n-1}y_n, x_{n-1}y_{n+1}, x_ny_{n-3}, x_ny_{n+1}, x_{n+1}y_{n-2}, x_{n+1}y_{n-1}\} \subseteq E(\overline{G})\]
	Therefore one can check that $C_8 =x_2y_nx_{n-1}y_{n+1}x_ny_{n-3}x_{n+1}y_{n-2}x_2\subseteq \overline{G}$, a contradiction again. So, assume that $xy\in \{x_{n-2}y_{n-2},x_{n-1}y_{n-2},x_{n-1}y_{n-1}\}$, that is $xy\in E(C)$.
	Hence the proof is the same as the case that $x=x_{n-1}, y=y_{n-1}$, where $xy\in E(C')$, now we can check that $C_{2n}\subseteq G$, which is a contradiction.
	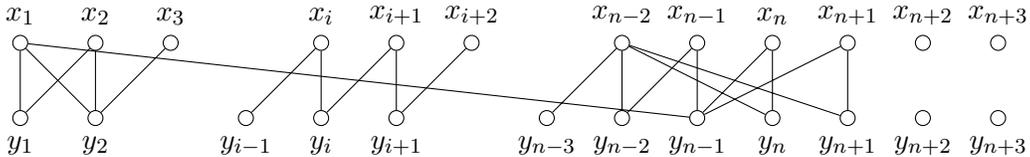
\begin{figure}[ht]
		\begin{center}
		\begin{tabular}{ccc}
			\begin{tikzpicture}
			\node [draw, circle, fill=white, inner sep=2pt, label=below:$y_1$] (y1) at (0,0) {};
			\node [draw, circle, fill=white, inner sep=2pt, label=below:$y_2$] (y2) at (1,0) {};
			\node [draw, circle, fill=white, inner sep=2pt, label=below:$y_{i-1}$](y3) at (3,0) {};
			\node [draw, circle, fill=white, inner sep=2pt, label=below:$y_{i} $] (y4) at (4,0) {};
			\node [draw, circle, fill=white, inner sep=2pt, label=below:$ y_{i+1}$] (y5) at (5,0) {};
			\node [draw, circle, fill=white, inner sep=2pt, label=below:$ y_{n-3}$](y6) at (7,0) {};
			\node [draw, circle, fill=white, inner sep=2pt, label=below:$y_{n-2}$] (y7) at (8,0) {};
			\node [draw, circle, fill=white, inner sep=2pt, label=below:$y_{n-1}$] (y8) at (9,0) {};
			\node [draw, circle, fill=white, inner sep=2pt, label=below:$y_{n}$] (y9) at (10,0) {};
			\node [draw, circle, fill=white, inner sep=2pt, label=below:$y_{n+1}$] (y10) at (11,0) {};
			\node [draw, circle, fill=white, inner sep=2pt, label=below:$y_{n+2}$] (y11) at (12,0) {};
			\node [draw, circle, fill=white, inner sep=2pt, label=below:$y_{n+3}$] (y12) at (13,0) {};
			\
			\node [draw, circle, fill=white, inner sep=2pt, label=above:$x_1$] (x1) at (0,1) {};
			\node [draw, circle, fill=white, inner sep=2pt, label=above:$x_2$] (x2) at (1,1) {};
			\node [draw, circle, fill=white, inner sep=2pt, label=above:$x_3$] (x3) at (2,1) {};
			\node [draw, circle, fill=white, inner sep=2pt, label=above:$ x_{i}$] (x4) at (4,1) {};
			\node [draw, circle, fill=white, inner sep=2pt, label=above:$ x_{i+1}$] (x5) at (5,1) {};
			\node [draw, circle, fill=white, inner sep=2pt, label=above:$ x_{i+2}$] (x6) at (6,1) {};
			\node [draw, circle, fill=white, inner sep=2pt, label=above:$x_{n-2}$] (x7) at (8,1) {};
			\node [draw, circle, fill=white, inner sep=2pt, label=above:$x_{n-1}$] (x8) at (9,1) {};
			\node [draw, circle, fill=white, inner sep=2pt, label=above:$x_n$] (x9) at (10,1) {};
			\node [draw, circle, fill=white, inner sep=2pt, label=above:$x_{n+1}$] (x10) at (11,1) {};
			\node [draw, circle, fill=white, inner sep=2pt, label=above:$x_{n+2}$] (x11) at (12,1) {};
			\node [draw, circle, fill=white, inner sep=2pt, label=above:$x_{n+3}$] (x12) at (13,1) {};
			\draw (x1)--(y1)--(x2)--(y2)--(x3);
			\draw (x1)--(y2);
			\draw (y6)--(x7)--(y7)--(x8)--(y8)--(x1);
			\draw (y3)--(x4)--(y4)--(x5)--(y5)--(x6);
			\draw (x7)--(y9)--(x9)--(y8)--(x10)--(y10)--(x7);
			\end{tikzpicture}\\
		\end{tabular}\\
		\caption{$(|X_1'|,|Y_1'|)=(1,1)$ and $x_{n-2}y_{n-1}\notin E(C')$}
		\label{fi7}
		\end{center}
	\end{figure}
	
	\bigskip
	Subcase 1-4. $(|X_1'|,|Y_1'|)=(2,0)$.\\
	Assume that $V(C') =\{x_{n-2},x_{n-1},x_n, y_n,y_{n+1},y_{n+2}\}$ and $C' =x_{n-2}y_nx_{n-1}y_{n+1}x_ny_{n+2}x_{n-2}$. In this case, we have:
	\[C_{2n} = x_1y_1\ldots y_{n-3}x_{n-2}y_{n+2}x_ny_{n+1}x_{n-1}y_{n-1}x_1\subseteq G\]
	Which is a contradiction.
	
	\bigskip
	Subcase 1-5. $(|X_1'|,|Y_1'|)=(2,1)$.\\
	W.l.g., we may assume that $V(C') =\{x_{n-2},x_{n-1},x_n, y,y_{n},y_{n+1}\}$, where $y\in \{y_{n-2},y_{n-1}\}$. In this case, by considering the edges of $C'$ and using Lemma \ref{le1},  it is easy to check that in any case $C_{2n}\subseteq G$, unless for the case that $y=y_{n-1}$ and $x_{n-2}y_{n-1}\in E(\overline{G})$. W.l.g., we may assume that $C'=x_{n-1}y_{n-1}x_ny_nx_{n-2}y_{n+1}x_{n-1}$. Consider the Figure \ref{fi8}.\\
	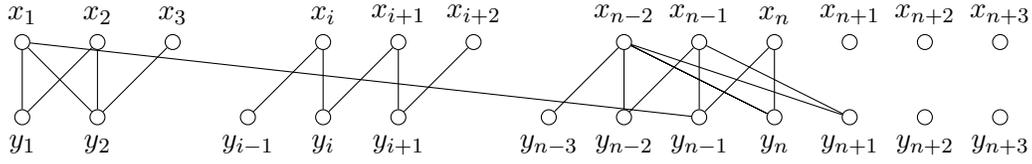
\begin{figure}[ht]
		\begin{center}
		\begin{tabular}{ccc}
			\begin{tikzpicture}
			\node [draw, circle, fill=white, inner sep=2pt, label=below:$y_1$] (y1) at (0,0) {};
			\node [draw, circle, fill=white, inner sep=2pt, label=below:$y_2$] (y2) at (1,0) {};
			\node [draw, circle, fill=white, inner sep=2pt, label=below:$y_{i-1}$](y3) at (3,0) {};
			\node [draw, circle, fill=white, inner sep=2pt, label=below:$y_{i} $] (y4) at (4,0) {};
			\node [draw, circle, fill=white, inner sep=2pt, label=below:$ y_{i+1}$] (y5) at (5,0) {};
			\node [draw, circle, fill=white, inner sep=2pt, label=below:$ y_{n-3}$](y6) at (7,0) {};
			\node [draw, circle, fill=white, inner sep=2pt, label=below:$y_{n-2}$] (y7) at (8,0) {};
			\node [draw, circle, fill=white, inner sep=2pt, label=below:$y_{n-1}$] (y8) at (9,0) {};
			\node [draw, circle, fill=white, inner sep=2pt, label=below:$y_{n}$] (y9) at (10,0) {};
			\node [draw, circle, fill=white, inner sep=2pt, label=below:$y_{n+1}$] (y10) at (11,0) {};
			\node [draw, circle, fill=white, inner sep=2pt, label=below:$y_{n+2}$] (y11) at (12,0) {};
			\node [draw, circle, fill=white, inner sep=2pt, label=below:$y_{n+3}$] (y12) at (13,0) {};
			\
			\node [draw, circle, fill=white, inner sep=2pt, label=above:$x_1$] (x1) at (0,1) {};
			\node [draw, circle, fill=white, inner sep=2pt, label=above:$x_2$] (x2) at (1,1) {};
			\node [draw, circle, fill=white, inner sep=2pt, label=above:$x_3$] (x3) at (2,1) {};
			\node [draw, circle, fill=white, inner sep=2pt, label=above:$ x_{i}$] (x4) at (4,1) {};
			\node [draw, circle, fill=white, inner sep=2pt, label=above:$ x_{i+1}$] (x5) at (5,1) {};
			\node [draw, circle, fill=white, inner sep=2pt, label=above:$ x_{i+2}$] (x6) at (6,1) {};
			\node [draw, circle, fill=white, inner sep=2pt, label=above:$x_{n-2}$] (x7) at (8,1) {};
			\node [draw, circle, fill=white, inner sep=2pt, label=above:$x_{n-1}$] (x8) at (9,1) {};
			\node [draw, circle, fill=white, inner sep=2pt, label=above:$x_n$] (x9) at (10,1) {};
			\node [draw, circle, fill=white, inner sep=2pt, label=above:$x_{n+1}$] (x10) at (11,1) {};
			\node [draw, circle, fill=white, inner sep=2pt, label=above:$x_{n+2}$] (x11) at (12,1) {};
			\node [draw, circle, fill=white, inner sep=2pt, label=above:$x_{n+3}$] (x12) at (13,1) {};
			\draw (x1)--(y1)--(x2)--(y2)--(x3);
			\draw (x1)--(y2);
			\draw (y6)--(x7)--(y7)--(x8)--(y8)--(x1);
			\draw (y3)--(x4)--(y4)--(x5)--(y5)--(x6);
			\draw (x7)--(y9)--(x7)--(y10)--(x8);
			\draw (y8)--(x9)--(y9);
			\end{tikzpicture}\\
		\end{tabular}\\
		\caption{$(|X_1'|,|Y_1'|)=(2,1)$, $y=y_{n-1}$ and $x_{n-2}y_{n-1}\in E(\overline{G})$}
		\label{fi8}
		\end{center}
	\end{figure}

	By Figure \ref{fi8}, it is easy to check that:
	\[\{x_1y_{n-2},x_1y_n,x_1y_{n+1},x_2y_{n-2},x_2y_n,x_2y_{n+1},x_{n-1}y_n, x_ny_{n-3}, x_ny_{n-2}, x_ny_{n+1} \}\subseteq E(\overline{G})\]
	Otherwise, we have $C_{2n}\subseteq G$. Now consider $x_{n-1}y_{n-3}$. If $x_{n-1}y_{n-3}\in E(\overline{G})$, we have $C_8=x_1y_{n-2}x_ny_{n-3}x_{n-1}y_nx_2y_{n+1}x_1\subseteq \overline{G}$, a contradiction. Hence assume that $x_{n-1}y_{n-3}\in E(G)$, therefore we have $C_{2n}=x_1y_{n-1}x_ny_nx_{n-2}y_{n-2}x_{n-1}y_{n-3}x_{n-3}y_{n-4} \ldots x_3y_2x_2y_1x_1\subseteq G$, which is a contradiction again.
	
	\bigskip
	Subcase 1-6. $(|X_1'|,|Y_1'|)=(2,2)$.\\
	W.l.g., we may assume that $V(C') =\{x_{n-2},x_{n-1},x_n, y_{n-2},y_{n-1},y_{n}\}$. If $x_ny_n\notin E(C')$, then $x_ny_j, y_nx_j \in E(C')$ for $j\in\{n-2,n-1\}$ and by Lemma \ref{le1}, we have $C_{2n}\subseteq G$. Hence assume that $x_ny_n\in E(C')$. Hence, it can be concluded that $|N_G(x_n)\cap \{y_{n-2}, y_{n-1}\}|\geq 1$ and $|N_G(y_n)\cap \{x_{n-2}, x_{n-1}\}|\geq 1$. Considering that $x_{n-2}y_{n-2}, y_{n-2}x_{n-1}, x_{n-1}y_{n-1}\in E(C)$,  if we have either $x_{n-1}y_n \in E(C')$ or $ x_{n}y_{n-2}\in E(C')$, then the proof follows the same logic, as indicated by Lemma \ref{le1}. Therefore, it can be say that $C' =x_{n-2}y_{n-2}x_{n-1}y_{n-1}x_ny_nx_{n-2} \subseteq G$. Since $C_{2n}\nsubseteq G$, by Figure \ref{fi9}, one can check that:
	\[\{x_1y_{n-2}, x_1y_n, x_2y_{n-2}, x_2y_n, x_{n-1}y_{n-3}, x_{n-1}y_n, x_ny_{n-3}, x_ny_{n-2}\}\subseteq E(\overline{G})\]
	
	Now we have the following claim:
	
	\bigskip
	\noindent\textbf{Claim 6.} $|N_G(x_n)\cap \{y_{n+1},y_{n+2},y_{n+3}\}|=0$. \\
	\begin{proof} By contradiction assume that $|N_G(x_n)\cap \{y_{n+1},y_{n+2},y_{n+3}\}|\neq0$, say $x_ny_{n+1}\in E(G)$. Hence $\{x_1,x_{n-1}\}\subseteq N_{\overline{G}}(y_{n+1})$; otherwise, $C_{2n}\subseteq G[X_1\cup \{x_n\}, Y_1\cup\{y_{n+1}\}]$ which is a contradiction. Since $\{x_1,x_{n-1}\}\subseteq N_{\overline{G}}(y_{n+1})$, we have $C_8 =x_1y_{n}x_2y_{n-2}x_ny_{n-3}x_{n-1}y_{n+1}x_1\subseteq \overline{G}$, a contradiction again. So the assumption does not hold and the claim holds.\end{proof}\\
	Therefore by Claim $6$, $\{y_{n+1},y_{n+2},y_{n+3}\}\subseteq N_{\overline{G}}(x_n)$. Now we have the following claim:
	
	\bigskip
	\noindent\textbf{Claim 7.} $\{x_1,x_2\}\subseteq N_{G}(y_i) $ for $i\in\{n+1,n+2,n+3\}$. \\
	\begin{proof} By contradiction we may assume that $x_1y_{n+1}\in E(\overline{G})$ (for other cases, the proof is identical). Therefore $C_8 =y_{n+1}x_1y_{n-2}x_2y_nx_{n-1}y_{n-3}x_ny_{n+1}\subseteq \overline{G}$, a contradiction.\end{proof}\\
	Now by Lemma \ref{le1} and  Claim $7$, for $i\in\{n+1,n+2,n+3\}$ we have $\{x_{n-1}y_i, x_{n-1}y_1,x_ny_1\}\subseteq E(\overline{G})$; otherwise, $C_{2n}\subseteq G$, a contradiction. If there exists a vertex $x$ of $\{x_{n+1},x_{n+2},x_{n+3}\}$, such that $|N_{\overline{G}}(x)\cap \{y_{n+1},y_{n+2},y_{n+3}\}|\geq 2$, then $C_8\subseteq \overline{G}$, a contradiction too. Hence, for each $i\in\{n+1,n+2,n+3\}$, we have $|N_G(x)\cap \{y_{n+1},y_{n+2},y_{n+3}\}|\geq 2$. W.l.g., we may assume that $x_{n+1}y_n, x_{n+1}y_{n+1}\in E(G)$ and thus one can check that $C_{2n}\subseteq G[X_1\cup\{x_{n+1}\}, (Y_1\setminus \{y_1\})\cup \{y_{n},y_{n+1}\}]$, which is a contradiction.\\
	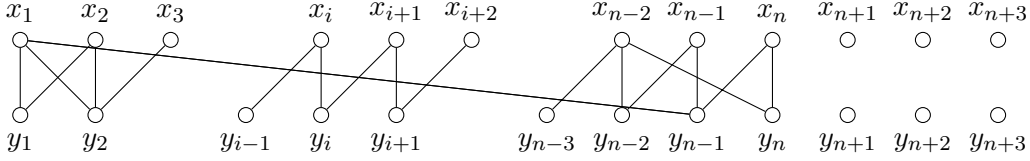
\begin{figure}[ht]
		\begin{center}
		\begin{tabular}{ccc}
			\begin{tikzpicture}
			\node [draw, circle, fill=white, inner sep=2pt, label=below:$y_1$] (y1) at (0,0) {};
			\node [draw, circle, fill=white, inner sep=2pt, label=below:$y_2$] (y2) at (1,0) {};
			\node [draw, circle, fill=white, inner sep=2pt, label=below:$y_{i-1}$](y3) at (3,0) {};
			\node [draw, circle, fill=white, inner sep=2pt, label=below:$y_{i} $] (y4) at (4,0) {};
			\node [draw, circle, fill=white, inner sep=2pt, label=below:$ y_{i+1}$] (y5) at (5,0) {};
			\node [draw, circle, fill=white, inner sep=2pt, label=below:$ y_{n-3}$](y6) at (7,0) {};
			\node [draw, circle, fill=white, inner sep=2pt, label=below:$y_{n-2}$] (y7) at (8,0) {};
			\node [draw, circle, fill=white, inner sep=2pt, label=below:$y_{n-1}$] (y8) at (9,0) {};
			\node [draw, circle, fill=white, inner sep=2pt, label=below:$y_{n}$] (y9) at (10,0) {};
			\node [draw, circle, fill=white, inner sep=2pt, label=below:$y_{n+1}$] (y10) at (11,0) {};
			\node [draw, circle, fill=white, inner sep=2pt, label=below:$y_{n+2}$] (y11) at (12,0) {};
			\node [draw, circle, fill=white, inner sep=2pt, label=below:$y_{n+3}$] (y12) at (13,0) {};
			\
			\node [draw, circle, fill=white, inner sep=2pt, label=above:$x_1$] (x1) at (0,1) {};
			\node [draw, circle, fill=white, inner sep=2pt, label=above:$x_2$] (x2) at (1,1) {};
			\node [draw, circle, fill=white, inner sep=2pt, label=above:$x_3$] (x3) at (2,1) {};
			\node [draw, circle, fill=white, inner sep=2pt, label=above:$ x_{i}$] (x4) at (4,1) {};
			\node [draw, circle, fill=white, inner sep=2pt, label=above:$ x_{i+1}$] (x5) at (5,1) {};
			\node [draw, circle, fill=white, inner sep=2pt, label=above:$ x_{i+2}$] (x6) at (6,1) {};
			\node [draw, circle, fill=white, inner sep=2pt, label=above:$x_{n-2}$] (x7) at (8,1) {};
			\node [draw, circle, fill=white, inner sep=2pt, label=above:$x_{n-1}$] (x8) at (9,1) {};
			\node [draw, circle, fill=white, inner sep=2pt, label=above:$x_n$] (x9) at (10,1) {};
			\node [draw, circle, fill=white, inner sep=2pt, label=above:$x_{n+1}$] (x10) at (11,1) {};
			\node [draw, circle, fill=white, inner sep=2pt, label=above:$x_{n+2}$] (x11) at (12,1) {};
			\node [draw, circle, fill=white, inner sep=2pt, label=above:$x_{n+3}$] (x12) at (13,1) {};
			\draw (x1)--(y1)--(x2)--(y2)--(x3);
			\draw (x1)--(y2);
			\draw (y6)--(x7)--(y7)--(x8)--(y8)--(x1);
			\draw (y3)--(x4)--(y4)--(x5)--(y5)--(x6);
			\draw (x7)--(y9)--(x9)--(y8)--(x1);
			\end{tikzpicture}\\
		\end{tabular}\\
		\caption{$(|X_1'|,|Y_1'|)=(2,2)$, $x_ny_n\in E(C')$}
		\label{fi9}
		\end{center}
	\end{figure}

	Case 2. For any $i\in\{1,2,\ldots,n-1\}$, each $x_iy_k \in E(\overline {G})$ if $k-i=1(\mod n-1)$ or $i-k=2(\mod n-1)$.\\
	Set $X_2 = \{x_{n-2},x_{n-1}, \ldots,x_{n+3}\}$ and $Y_2=\{y_{n-2},y_{n-1}, \ldots,y_{n+3}\}$. Since $BR(C_6, C_{8})=6$, $|X_2|=|Y_2|=6$, and $C_{8}\nsubseteq{\overline G}$, $G[X_2,Y_2]$ has a subgraph $C'\cong C_6$. Let $X'_1 =V(C')\cap X_1$ and $Y'_1 =V(C')\cap Y_1$. As the same as case 1,  we have the following subcases:
	
	\bigskip
	Subcase 2-1. $(|X_1'|,|Y_1'|)=(0,0)$.\\
	W.l.g., assume that $V(C') =\{x_{n},x_{n+1},x_{n+2}, y_{n},y_{n+1},y_{n+2}\}$ and $C' =x_{n}y_{n}x_{n+1}y_{n+1}x_{n+2}y_{n+2}x_{n}$. Set $X''_2 =\{x_{n},x_{n+1},x_{n+2}\}$, $Y''_2 =\{y_{n},y_{n+1},y_{n+2}\}$. Since $C_{8}\nsubseteq \overline{G}$, there is at least one edge between $X_1$ and $Y''_2$, or  at least one edge between $X''_2$ and $Y_1$. W.l.g., assume that $x_ny_{n-1}\in E(G)$. For other cases, the proof is the same. Hence we have $K_{3,2}\cong [\{x_{n-3}, x_{n-2}, x_{n-1}\}, \{y_n, y_{n+1}\}]\subseteq \overline{G}$. Similarly we have $y_1x_{n+1}, y_1x_{n+2}, y_{n-2}x_{n+1}, y_{n-2}x_{n+2} \in E(\overline{G})$. Therefore  $C_8 \subseteq \overline{G}$, a contradiction.
	
	
	\bigskip
	Subcase 2-2. $(|X_1'|,|Y_1'|)=(1,0)$.\\
	Assume that $V(C') =\{x_{n-1},x_n,x_{n+1}, y_{n},y_{n+1},y_{n+2}\}$ and $C' =x_{n-1}y_{n}x_{n}y_{n+1}x_{n+1}y_{n+2}x_{n-1}\subseteq G$. Now one can check that $K_{2,4}\cong[\{x_{n}, x_{n+1}\},\{y_1, y_{n-3}, y_{n-2}, y_{n-1}\}]\subseteq \overline{G}$ and $x_{n-3}y_n, x_{n-3}y_{n+2} \in E(\overline{G})$. Otherwise, we have $C_{2n}\subseteq G$, a contradiction. For example, by contrary assume that $x_ny_1\in E(G)$, so it can be seen that $C_{2n}:=x_ny_1x_2y_2x_2\ldots x_{n-2}y_{n-2}x_{n-1}y_{n+2}x_{n+1}y_{n+1}x_n\subseteq G$. For other cases, the proof is the same.\\
	Since $x_{n-4}y_{n-3}, x_{n-3}y_{n-2},x_{n-2}y_{n-1}\in E(\overline{G})$, if $x_{n-2}y_n$ or $x_{n-2}y_{n+2}$ belong to $E(\overline{G})$, then we have $C_8\subseteq \overline{G}$. Hence assume that $x_{n-2}y_n, x_{n-2}y_{n+2}\in E(G)$ and thus $x_{n-4}y_n\in E(\overline{G})$, if not, $C_{2n}\subseteq G$. So $C_8=y_1x_{n+1}y_{n-2}x_{n-3}y_nx_{n-4}y_{n-3}x_ny_1\subseteq \overline{G}$, a contradiction too.
	
	\bigskip
	Subcase 2-3. $(|X_1'|,|Y_1'|)=(2,0)$.\\
	Let us make the assumption, without loss of generality, that $V(C') =\{x_{n-2},x_{n-1},x_n, y_n,y_{n+1},y_{n+2}\}$. Furthermore, leveraging  of symmetry and without loss of generality, we can assume that $C'$ has a form like $C' =x_{n-2}y_nx_{n-1}y_{n+1}x_ny_{n+2}x_{n-2}$.  In this scenario, the following contradiction arises:
	\[C_{2n} = x_1y_1x_2y_2\ldots x_{n-3}y_{n-3}x_{n-2}y_{n+2}x_ny_{n+1}x_{n-1}y_{n-1}x_1\subseteq G\]
	
	\bigskip
	Subcase 2-4. $(|X_1'|,|Y_1'|)=(1,1)$.\\
	W.l.g., assume that $V(C') =\{x,x_{n},x_{n+1}, y, y_{n},y_{n+1}\}$, where $x\in \{x_{n-2},x_{n-1}\}$ and $ y\in \{y_{n-2},y_{n-1}\} $. If $xy\notin E(C)$, then $x=x_{n-2}$ and $y=y_{n-1}$. Let $x_{n-2}y_{n-1}\in E(C')$ and assume that $C' =x_{n-2}y_{n-1}x_ny_nx_{n+1}y_{n+1}x_{n-2}$. In this case, we have $C_{2n} = x_1y_1x_2y_2\ldots x_{n-2}y_{n+1}x_{n+1}y_nx_ny_{n-1}x_1\subseteq G$, a contradiction. Assume that $x_{n-2}y_{n-1}\notin E(C')$ and w.l.g. let $C' =x_{n-2}y_nx_ny_{n-1}x_{n+1}y_{n+1}x_{n-2}$. Since $C_{2n}  \nsubseteq G$, it is easy to check that: 
	\[\{ x_1y_n, x_1y_{n+1}, x_{n-1}y_n, x_{n-1}y_{n+1}, x_ny_{n-2}, x_ny_{n-3}, x_{n+1}y_{n-2}, x_{n+1}y_{n}, x_{n}y_{n+1} \} \subseteq E(\overline{G})\]
	
	For example, by contrary assume that $x_1y_n\in E(G)$, hence it can be checked that $G$ has at least one copy of $C_{2n}$, namely 
	$C_{2n}:=x_1y_1x_2y_2\ldots x_{n-3}y_{n-3}x_{n-2}y_{n+1}x_{n+1}y_{n-1}x_ny_nx_1$. For other cases, the proof is the same.\\
	Since $x_{n-1}y_{n-3}\in E(\overline{G})$, we have $C_8 =x_1y_{n+1}x_{n-1}y_{n-3}x_ny_{n-2}x_{n+1}y_{n}x_1\subseteq \overline{G}$, a contradiction again. Therefore, assume that $xy\in \{x_{n-2}y_{n-2},x_{n-1}y_{n-2},x_{n-1}y_{n-1}\}$, that is $xy\in E(C)$.
	Hence the proof is the same as the case that $x=x_{n-1}, y=y_{n-1}$ where $xy\in E(C')$,that is we can check that $C_{2n}\subseteq G$, a contradiction.

	\bigskip
	Subcase 2-5. $(|X_1'|,|Y_1'|)=(2,1)$.\\
	W.l.g., we may assume that $V(C') =\{x_{n-2},x_{n-1},x_n, y,y_{n},y_{n+1}\}$ where $y\in \{y_{n-2},y_{n-1}\}$. In this case by considering the edges of $C'$ and by Lemma \ref{le1},  it is easy to check that $C_{2n}\subseteq G$ for any case, unless for the case that $y=y_{n-1}$ and $x_{n-2}y_{n-1}\in E(\overline{G})$. W.l.g., we may assume that $C'=x_{n-1}y_{n-1}x_ny_nx_{n-2}y_{n+1}x_{n-1}$.  Now it is easy to check that:
	\[\{x_1y_{n-2},x_1y_n,x_1y_{n+1},x_2y_n,x_2y_{n+1},x_{n-1}y_n, x_ny_{n-3}, x_ny_{n-2}, x_ny_{n+1} \}\subseteq E(\overline{G})\]
	Otherwise,  $C_{2n}\subseteq G$. For example, by contrary assume that $x_1y_{n-2}\in E(G)$, hence it can be checked that $G$ has at least one copy of $C_{2n}$, namely $C_{2n}:=x_1y_1x_2y_2\ldots x_{n-3}y_{n-3}x_{n-2}y_{n}x_{n}y_{n-1}x_{n-1}y_{n-2}x_1$. For other cases, the proof is the same.\\
	Since $x_{n-1}y_{n-3}\in E(\overline{G})$, we have $C_8=x_1y_{n-2}x_ny_{n-3}x_{n-1}y_nx_2y_{n+1}x_1\subseteq \overline{G}$, a contradiction again.

	\bigskip 
	Subcase 2-6. $(|X_1'|,|Y_1'|)=(2,2)$.\\
	W.l.g., we may assume that $V(C') =\{x_{n-2},x_{n-1},x_n, y_{n-2},y_{n-1},y_{n}\}$. If $x_ny_n\notin E(C')$, then $x_ny_j, y_nx_j \in E(C')$ for $j\in\{n-2,n-1\}$, and by Lemma \ref{le1} we have $C_{2n}\subseteq G$. Hence assume that $x_ny_n\in E(C')$. Now, one can say that $|N_G(x_n)\cap \{y_{n-2}, y_{n-1}\}|\geq 1$ and $|N_G(y_n)\cap \{x_{n-2}, x_{n-1}\}|\geq 1$. Now, since $x_{n-2}y_{n-2}, y_{n-2}x_{n-1}, x_{n-1}y_{n-1}\in E(C)$, if $x_{n-1}y_n \in E(C')$ or $ x_{n}y_{n-2}\in E(C')$, then by Lemma \ref{le1}, the proof is complete. Therefore, one can assume that $C' =x_{n-2}y_{n-2}x_{n-1}y_{n-1}x_ny_nx_{n-2} \subseteq G$. Since $C_{2n}\nsubseteq G$, one can check that:
	\[\{x_1y_{n-2}, x_1y_n, x_{n-1}y_{n-3}, x_{n-1}y_n, x_ny_{n-3}, x_ny_{n-2}\}\subseteq E(\overline{G})\]
	Set $X_3=\{x_{n+1}, x_{n+2}, x_{n+3}\}$, $Y_3=\{y_{n+1}, y_{n+2}, y_{n+3}\}$. There exists at least one vertex of $X_3$ or $Y_3$, say $y_{n+1}$, such that $|N_{\overline{G}}(y_{n+1}) \cap X_3|\geq 2$; otherwise, we have $C_6\subseteq G[X_3,Y_3]$ and the proof is the same as subcase 2-1. Assume that $x_{n+1}y_{n+1}, x_{n+2}y_{n+1} \in E(\overline{G})$. Therefore, for each $x\in \{x_{n+1}, x_{n+2}\}$, we have $|N_{\overline{G}}(x) \cap \{y_{n-3}, y_{n-2}, y_n\}|\geq 2$; otherwise, $C_{2n}\subseteq G$, a contradiction. If $|N_{\overline{G}}(x) \cap \{y_{n-3}, y_{n-2}, y_n\}|\geq 2$ for some $x\in \{x_{n+1}, x_{n+2}\}$ or $N_{\overline{G}}(x_{n+1}) \cap \{y_{n-3}, y_{n-2}, y_n\}\neq N_{\overline{G}}(x_{n+2}) \cap \{y_{n-3}, y_{n-2}, y_n\}$, then we have $C_8\subseteq \overline{G}$, a contradiction. That is, we have $|N_{\overline{G}}(x) \cap \{y_{n-3}, y_{n-2}, y_n\}|= 1$ for each $x\in \{x_{n+1}, x_{n+2}\}$ and one can check that $N_{\overline{G}}(x) \cap \{y_{n-3}, y_{n-2}, y_n\}=\{y_n\}$. Therefore we have $x_1y_{n+1}, x_{n-1}y_{n+1}\in E(G)$, if not, $C_{2n}\subseteq G$, which is a contradiction. Now, we have $C_{n}=x_1y_1x_2y_2\ldots x_{n-}y_{n-3}x_{n-2}y_nx_ny_{n-1}x_{n-1}y_{n+1}x_1\subseteq G$, a contradiction again.
	
	Hence by Cases $1,2$, the proof is complete and the theorem holds.
\end{proof}

Therefore by Lemmas \ref{le3} and \ref{le4} and by Theorems \ref{th4}, \ref{th5}, \ref{th6} and \ref{th7}, Theorem \ref{th1} holds.\\


\nocite{*}
\bibliographystyle{abbrvnat}
\bibliography{sample-dmtcs}
\label{sec:biblio}

\end{document}